\documentclass[11pt,a4paper]{article}
\usepackage[]{times}
\usepackage{fullpage}
\usepackage{graphicx}
\usepackage{subfigure}
\usepackage{amsmath}
\usepackage{fancyhdr}
\date{}

\newcommand{\prov}{{\sc Proof}.\hspace*{3mm} }

\newcommand{\QED}{$\rule{2mm}{2mm}$}

\newtheorem{theorem}{Theorem}[section]
\newtheorem{lemma}[theorem]{Lemma}
\newtheorem{e-proposition}[theorem]{Proposition}
\newtheorem{corollary}[theorem]{Corollary}
\newtheorem{e-definition}[theorem]{Definition\rm}
\newtheorem{remark}{\it Remark\/}

\title{Optimal simplex finite-element approximations of arbitrary order in curved domains circumventing the isoparametric technique}

\author{
    Vitoriano Ruas$^{1,2}$\thanks{This work was partially supported by CNPq, the National Research Council of Brazil}
		\\[1mm]
  {\small $^{1}$ Sorbonne Universit\'e, UMR 7190, Institut Jean Le Rond d'Alembert, France}\\
  {\small $^{2}$ CNRS, UMR 7190, Institut Jean Le Rond d'Alembert, F-75005 Paris, France.}\\[1mm]
  {\small e-mail: {\it vitoriano.ruas@upmc.fr}}}
	
\begin{document}
\maketitle
\thispagestyle{fancy}

\begin{abstract}

Since the 1960's the finite element method emerged as a powerful tool for the numerical simulation of countless physical phenomena or processes in 
applied sciences. One of the reasons for this undeniable success is the great versatility of the finite-element approach to deal with different types of 
geometries. This is particularly true of problems posed in curved domains of arbitrary shape. In the case of function-value Dirichlet conditions prescribed on curvilinear boundaries method's isoparametric version for meshes 
consisting of curved triangles or tetrahedra has been mostly employed to recover the optimal approximation properties known to hold for standard 
straight elements in the case of polygonal or polyhedral domains. However, besides obvious algebraic and geometric inconveniences, the isoparametric technique 
is helplessly limited in scope and simplicity, since its extension to degrees of freedom other than function values is not straightforward if not unknown. 
The purpose of this paper is to propose, study and test a simple alternative that bypasses all the above drawbacks, without eroding qualitative approximation
properties. More specifically this technique can do without curved elements and is based only on polynomial algebra. 
\end{abstract}

\section{Study framework}

This work deals with a new method for solving boundary value problem posed in a two- or three-dimensional domain, with a smooth curved boundary of arbitrary shape. In the framework of the finite-element solution of second order elliptic equations posed in curved domains with Dirichlet boundary conditions, it is well known that 
a considerable order lowering may occur if prescribed boundary values are shifted to nodes that are not mesh vertexes of an approximating polygon or polyhedron 
formed by the union of the ordinary $N$-simplexes of a fitted mesh. Over four decades ago some techniques were designed in order to remedy such a loss of accuracy, and possibly attain the same theoretical optimal orders as in the case of a polytopic domain, assuming that the solution is sufficiently smooth. Two examples of such attempts are the \textit{interpolated boundary condition method} by Nitsche and Scott (cf. \cite{Nitsche} and \cite{Scott}), and the method introduced by Zl\'amal in \cite{Zlamal} and extended by \v{Z}\'eni\v{s}ek in \cite{Zenisek}.\\
The principle our method is based upon is close to the interpolated boundary conditions 
studied in \cite{BrennerScott} for two-dimensional problems. Although the latter technique is very intuitive and has been known since the seventies (cf. 
\cite{Scott}), it has been of limited use so far. Among the reasons for this we could quote its difficult implementation, the lack of an extension to three-dimensional problems, and most of all, restrictions on the choice of boundary nodal points to reach optimal convergence rates. In contrast our method is simple to implement in both two- and three-dimensional geometries. Moreover optimality is attained very naturally in both cases for various choices of boundary nodal points. \\
\indent In order to allow an easier description of our methodology we consider as a model the convection-diffusion equation with Dirichlet boundary conditions, solved by different $N$-simplex based  methods, incorporating degrees of freedom other than function values at the mesh vertexes. For instance, 
if standard quadratic Lagrange finite elements are employed, it is well-known that approximations of an order not greater than $1.5$ in the energy norm are generated (cf. \cite{Ciarlet}), in contrast to the second order ones that apply to the case of a polygonal or polyhedral domain, assuming that the solution is sufficiently smooth. If we are to recover the optimal second order approximation property something different has to be done. \\
Since long the isoparametric version of the finite element method for meshes consisting of curved triangles or tetrahedra (cf. \cite{Zienkiewicz}), has been considered as the ideal way to achieve this. It turns out that, besides a more elaborated description of the mesh, the isoparametric technique inevitably leads to the integration of rational functions to compute the system 
matrix. This raises the delicate question on how to choose the right numerical quadrature formula in the master element, especially in the case of complex non linear problems. In contrast, in the technique to be introduced in this paper exact numerical integration can always be used for this purpose, since we only have to deal with polynomial integrands. Moreover the element geometry remains the same as in the case of polygonal or polyhedral domains. It is noteworthy that both advantages are conjugated with the fact that no erosion of qualitative approximation properties results from the application of our technique, as compared to the equivalent isoparametric one. 
We should also emphasize that this approach is particularly handy, whenever the finite element method under consideration has normal components or normal derivatives 
as degrees of freedom. Indeed in this case isoparametric analogs are either not so easy to define (see. e.g. \cite{RT}) or are simply unknown. \\
\indent An outline of the paper is as follows. In Section 2 we present our method to solve the model problem with Dirichlet boundary conditions in a smooth curved two-dimensional domain with conforming Lagrange finite elements based on meshes with straight triangles, in connection with the standard Galerkin formulation. Corresponding well-posedness results are demonstrated. In Section 3 we prove general error estimates for the method introduced in the previous section. Moreover $L^2$-error estimates are demonstrated in relevant cases, which to the best of author's knowledge are unprecedented for the class of problems considered in this work.
In Section 4 we assess the approximation properties of the method studied in the previous section by solving some two-dimensional test-problems with 
piecewise quadratic functions. 
We conclude in Section 5 with some comments on the methodology studied in this work. In particular we briefly show that the technique addressed in Sections 2 and 3 applies with no particular difficulty to the case of boundary value problems posed in curved three-dimensional domains (see also \cite{arXiv3D}).


\section{Method description}
The methodology to enforce Dirichlet boundary conditions on curvilinear boundaries considered in this work applies to many types of equations. However, in order to avoid non essential difficulties, we consider as a model the following convection-diffusion equation in an $N$-dimensional smooth domain $\Omega$ with boundary $\Gamma$, for $N=2$ or $N=3$, namely:
\begin{equation}
\label{Poisson}
\left\{
\begin{array}{l}
 -\nu \Delta u + {\bf b} \cdot {\bf grad} \; u = f \mbox{ in } \Omega \\
 u = d \mbox{ on } \Gamma,
\end{array}
\right.
\end{equation} 
\noindent where $\nu$ is the diffusion coefficient and ${\bf b} \in [L^{\infty}(\Omega)]^N$ is a given convective velocity assumed to be divergence free. 
$f$ and $d$ in turn are given functions defined in $\Omega$ and on $\Gamma$, having suitable regularity properties. We shall be dealing with approximation 
methods of order $k$ for $k > 1$ in the standard energy norm $\parallel {\bf grad} (\cdot) \parallel_{0}$, as long as $u \in H^{k+1}(\Omega)$, where 
$\parallel \cdot \parallel_{0}$ equals $[ \int_{\Omega} (\cdot)^2 ]^{1/2}$, i.e. it denotes the standard norm of $L^2(\Omega)$.  
Accordingly, we shall assume that $f \in H^{k-1}(\Omega)$ and $d \in H^{k+1/2}(\Gamma)$ (cf. \cite{Adams}). Although the method to be described below applies to any 
$d$, for the sake of simplicity henceforth we shall take $d \equiv 0$. In this case, for the assumed regularity of $u$ to hold, we require that both ${\bf b}$ and $\Gamma$ be sufficiently smooth and at least of the $C^{k-1}$-class. \\ 
In order to simplify the presentation here we confine the description of our method to the two-dimensional case, leaving an overview of the 
three-dimensional case for Section 5.\\
\indent Let us be given a mesh ${\mathcal T}_h$ conssting of straight-edged triangles satisfying the usual compatibility conditions and fitting $\Omega$ 
in such a way that all the vertexes of the polygon $\Omega_h := \cup_{T \in {\mathcal T}_h}$ belong to $\Gamma$. 
Every element of ${\mathcal T}_h$ is considered to be a closed set and is assumed to belong to a uniformly regular family of partitions (see e.g. 
\cite{Ciarlet}). Let $\Gamma_h$ be the boundary of $\Omega_h$ and $h_T$ be the diameter of $T \in {\mathcal T}_h$. As usual we 
set $h :=  \max_{T \in {\mathcal T}_h} h_T$. Clearly enough if $\Omega$ is convex $\Omega_h$ is a proper subset of $\Omega$. We make the very reasonable assumptions on the mesh that no element in ${\mathcal T}_h$ has more than one edge on $\Gamma_h$. \\
\noindent We also need some definitions regarding the skin $(\Omega \setminus \Omega_h) \cup (\Omega_h \setminus \Omega)$. First of all, in order to avoid non essential technicalities, we 
assume that the mesh is constructed in such a way that convex and concave portions of $\Gamma$ correspond to convex and concave portions of $\Gamma_h$. This 
property is guaranteed if the points separating such portions of $\Gamma$ are vertexes of polygon $\Omega_h$. In doing so, let 
${\mathcal S}_h$ be the subset of ${\mathcal T}_h$ consisting of triangles having one edge on $\Gamma_h$. Now $\forall T \in {\mathcal S}_h$ we denote by  
$\Delta_T$ the closed set delimited by $\Gamma$ and the edge $e_T$ of $T$ whose end-points belong to $\Gamma$ and set $T^{'} := T \cup \Delta_T$ 
if $\Delta_T$ is not a subset of $T$ and $T^{'}:= \overline{T \setminus \Delta_T}$ otherwise (see Figure 1).
\begin{figure}[h]
\label{fig1}
\begin{center}
\includegraphics[scale=0.5]{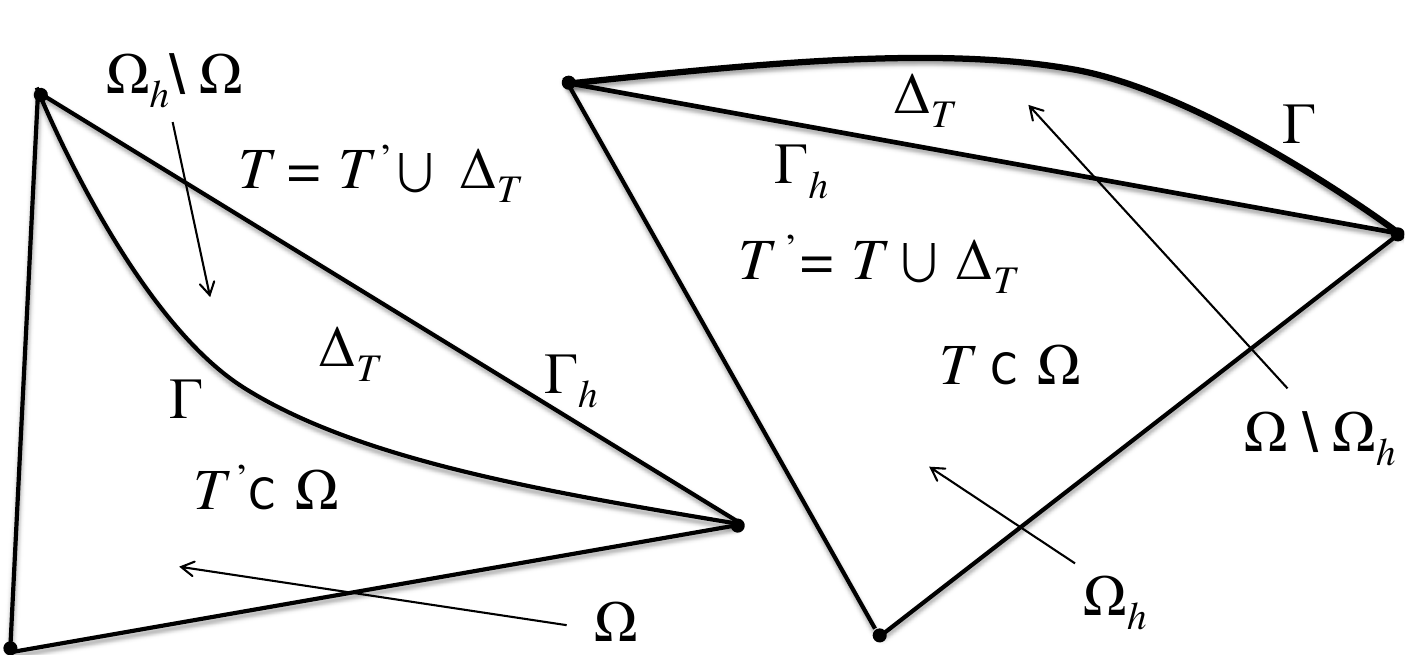}
\end{center}
\par
\caption{Skin $\Delta_T$ related to a mesh triangle $T$ next to a convex (right) or a concave (left) portion of $\Gamma$}
\end{figure} 
Notice that if $e_T$ lies on a convex portion of $\Gamma_h$, $T$ is a proper subset of $T^{'}$, while the opposite occurs if $e_T$ lies on a concave portion of $\Gamma_h$. With such a definition we can assert that there is a partition ${\mathcal T}_h^{'}$ of $\Omega$ associated with ${\mathcal T}_h$ consisting of non overlapping 
sets $T^{'}$ for $T \in {\mathcal S}_h$, besides the elements in ${\mathcal T}_h \setminus {\mathcal S}_h$. \\ 

\indent Next we introduce two spaces $V_h$ and $W_h$ associated with ${\mathcal T}_h$. $V_h$ is the standard Lagrange finite element space consisting of 
continuous functions $v$ defined in $\Omega_h$ that vanish on $\Gamma_h$, whose restriction to every $T \in {\mathcal T}_h$ is a polynomial of degree 
less than or equal to $k$ for $k \geq 2$. For convenience we extend by zero every function $v \in V_h$ to $\Omega \setminus \Omega_h$. $W_h$ in turn is the space 
of functions defined in $\Omega_h$ having the properties listed below. 
 
\begin{enumerate} 
\item The restriction of $w \in W_h$ to every $T \in {\mathcal T}_h$ is a polynomial of degree less than or equal to $k$;
\item Every $w \in W_h$ is continuous in $\Omega_h$ and vanishes at the vertexes of $\Gamma_h$; 
\item A function $w \in W_h$ is extended to $\Omega \setminus \Omega_h$ in such a way that its polynomial expression in $T \in {\mathcal S}_h$ also applies 
to points in $\Delta_T$;
\item $\forall T \in {\mathcal S}_h$, $w(P) = 0$ for every $P$ among the $k-1$ nearest intersections with $\Gamma$ 
of the line passing through the vertex $O_T$ of $T$ not belonging to $\Gamma$ and the points $M$ different from vertexes of $T$ subdividing the edge $e_T$ opposite 
to $O_T$ into $k$ segments of equal length (cf. Figure 2).
\end{enumerate}
\begin{figure}[h]
\label{fig2}
\begin{center}
\includegraphics[scale=0.5]{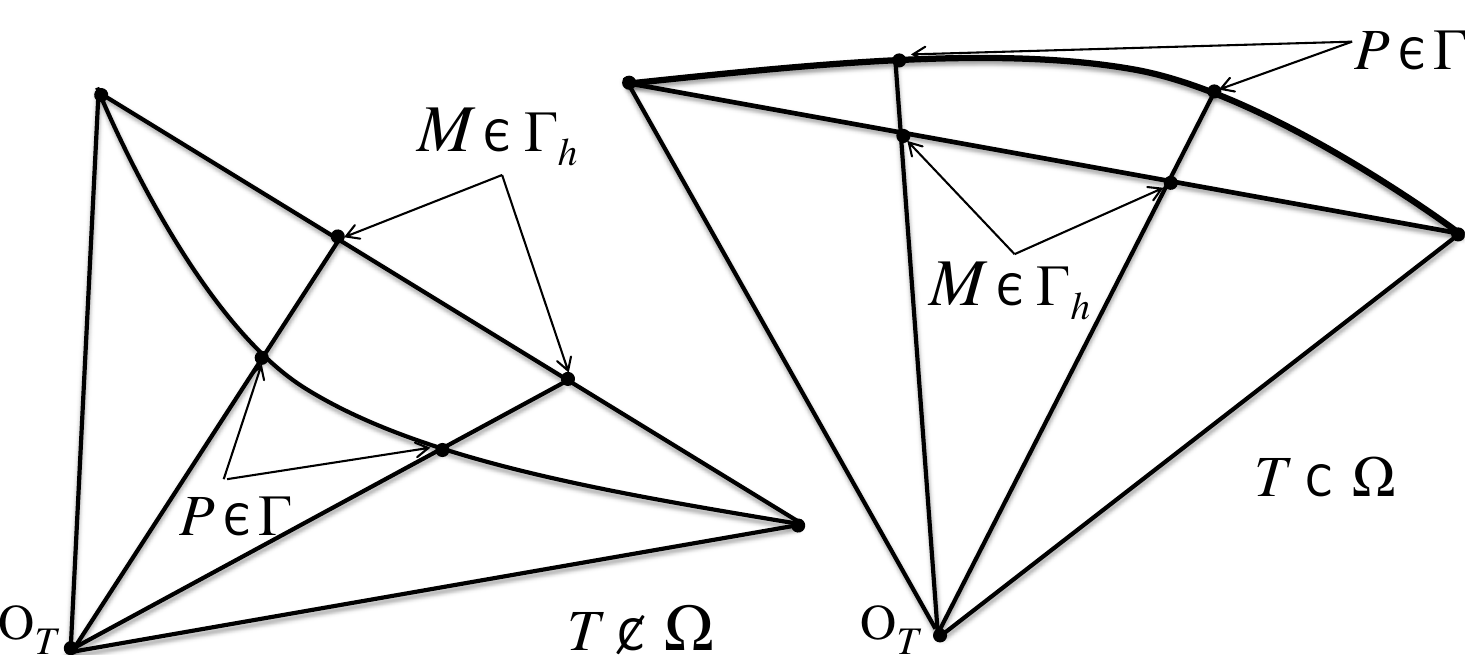}
\end{center}
\par
\caption{Construction of nodes $P \in \Gamma$ for space $W_h$ related to Lagrangian nodes $M \in \Gamma_h$ for $k=3$}
\end{figure}  
\begin{remark}
The construction of the nodes associated with $W_h$ located on $\Gamma$ advocated in item 4 is not mandatory. Notice that it differs 
from the intuitive construction of such nodes lying on normals to edges of $\Gamma_h$ commonly used in the isoparametric technique. The main advantage of this proposal is an easy determination of boundary node coordinates by linearity, using a supposedly available analytical expression of $\Gamma$. 
Nonetheless the choice of boundary nodes ensuring our method's optimality is really wide, in contrast to the restrictions inherent to 
the interpolated boundary condition method (cf. \cite{BrennerScott}).  \QED
\end{remark}
The fact that $W_h$ is a non empty finite-dimensional space is next established.

\begin{lemma}
\label{lemma1}
Let ${\mathcal P}_k(T)$ be the space of polynomials defined in $T \in {\mathcal S}_h$ of degree less than or equal to $k$.   
Provided $h$ is small enough $\forall T \in {\mathcal S}_h$, given a set of $m_k$ real values $b_{i}$, $i=1,\ldots,m_k$ with $m_k=(k+1)k/2$, 
there exists a unique function $w_T \in {\mathcal P}_k(T)$ that vanishes at both 
vertexes of $T$ located on $\Gamma$ and at the $k-1$ points $P$ of $\Gamma$ defined in accordance with item 4. of the above definition of $W_h$, and takes 
value $b_i$ respectively at the $m_k$ nodes of $T$ not located on $\Gamma_h$, corresponding to the Lagrange family of triangular finite elements  
(cf. \cite{Zienkiewicz}). 
\end{lemma}   
 
\prov Let us first extend the vector $\vec{b}:=[b_1,b_2,\ldots,b_{m_k}]$ of $\Re^{m_k}$ into a vector of $\Re^{n_k}$ still denoted by $\vec{b}$, with $n_k:=m_k+k+1$,  by adding $n_k-m_k$ zero components. If the boundary nodes $P$ were replaced by the corresponding $M \in \Gamma_h \cap T$, it is clear that the result would hold true, according to the well-known properties of Lagrange finite elements. The vector $\vec{a}$ of coefficients $a_i$ for $i=1,2,\ldots,n_k := (k+2)(k+1)/2$ of the canonical basis functions 
$\varphi_i$ of ${\mathcal P}_k(T)$ for $1 \leq i \leq n_k$ would be given by $a_i=b_i$ for $1 \leq i \leq n_k$. Denoting by $M_i$ the Lagrangian nodes of $T$, $i=1,2,\ldots,n_k$, this means that the matrix $K$ whose entries are $k_{ij} := \varphi_j(M_i)$ is the identity matrix. Let $\tilde{M}_i=M_i$ if $M_i \notin \Gamma \setminus \Gamma_h$ and $\tilde{M}_i$ be the node of the type $P$ associated with $M_i$ otherwise. The Lemma will be proved if the $n_k \times n_k$ linear system $\tilde{K} \vec{a} = \vec{b}$ is uniquely solvable, where $\tilde{K}$ is the matrix with entries $\tilde{k}_{ij}:=\varphi_j(\tilde{M}_{i})$. Clearly we have $\tilde{K} = K + 
E_K$, where the entries of $E_K$ are $e_{ij}:= \varphi_j(\tilde{M}_{i}) - \varphi_j(M_{i})$. At this point we observe that there exists a constant $C_{\Gamma}$ depending only on $\Omega$ such that the length of the segment $\overline{M_i\tilde{M}_i}$ is bounded above by $C_{\Gamma} h_T^2$. It follows that $\forall \;i,j$, 
$|e_{ij}| \leq C_{\Gamma} h_T ^2 \max_{{\bf x} \in T \cup \Delta_T} |{\bf grad}\;\varphi_j({\bf x})|$.  \\
Since $\varphi_j$ is a polynomial and $\Delta_T$ is at most a small perturbation of $T$, the maximum of $|{\bf grad}\;\varphi_j|$ in $T \cup \Delta_T$ must be bounded by a certain mesh independent constant times $\max_{{\bf x} \in T} |{\bf grad}\;\varphi_j({\bf x})|$. From standard arguments we know that the latter maximum is bounded above by a mesh-independent constant times $h_T^{-1}$. In short we have $|e_{ij}| \leq C_E h_T$ $\forall \; i,j$, where $C_E$ is a mesh independent constant. Hence the matrix $\tilde{K}$ equals the identity matrix plus an $O(h_T)$ matrix $E_K$. Therefore $\tilde{K}$ is an invertible matrix, as long as $h$ is sufficiently small. \QED \\         
  
Now let us set the problem associated with spaces $V_h$ and $W_h$, whose solution is an approximation of $u$, that is, the solution of (\ref{Poisson}). Extending 
$f$ in $\Omega_h \setminus \Omega$ in different ways to be specified hereafter, and still denoting the resulting function defined in $\Omega \cup \Omega_h$ by $f$, 
we wish to solve,
\begin{equation}
\label{Poissonh}
\left\{
\begin{array}{l}
\mbox{Find } u_h \in W_h \mbox{ such that} \\ 
a_h(u_h,v) = F_h(v) \; \forall v \in V_h \\
\mbox{where } a_h(w,v) := \int_{\Omega_h} [\nu {\bf grad}\; w \cdot {\bf grad}\; v + ({\bf b} \cdot {\bf grad}\;w) v] \mbox{ and } F_h(v) : = \int_{\Omega_h} f v.
\end{array}
\right.
\end{equation}  

For convenience henceforth we refer to the nodes in a triangle belonging to the set of $(k+2)(k+1)/2$ points used to define the space of polynomials of degree less than or equal to $k >1$ for Lagrange finite elements, as the \textit{Lagrangian nodes} 
(cf. \cite{Ciarlet}, \cite{Zienkiewicz}). \\
Let us denote by $\parallel \cdot \parallel_{0,h}$ the standard norm of $L^2(\Omega_h)$. We next prove: 
\begin{e-proposition}
\label{propo1}
Provided $h$ is sufficiently small problem (\ref{Poissonh}) has a unique solution. Moreover there exists a constant $\alpha > 0$ independent of $h$ such that,
\begin{equation}
\label{infsup}
\forall w \in W_h \neq 0, \displaystyle \sup_{v \in V_h \setminus \{ 0 \}} \frac{a_h(w,v)}{\parallel {\bf grad} \; w \parallel_{0,h} \parallel {\bf grad} \; v \parallel_{0,h}} 
\geq \alpha.
\end{equation}  
\end{e-proposition}

\prov Given $w \in W_h$ let $v \in V_h$ coincide with $w$ at all Lagrangian nodes of elements $T \in {\mathcal T}_h \setminus {\mathcal S}_h$. As for an element 
$S \in {\mathcal S}_h$ we set $v=w$ at the Lagrangian nodes not belonging to $\Gamma_h$ and $v=0$ at the Lagrangian nodes located on $\Gamma_h$. 
The fact that on the edges common to two mesh elements $T^{-}$ and $T^{+}$, both $v_{|T^{-}}$ and $v_{|T^{+}}$ are polynomials of degree less than or equal to $k$ 
in terms of one variable coinciding at the exact number of points required to uniquely define such a function, implies that $v$ is continuous in $\Omega_h$. Moreover for the same reason $v$ vanishes all over $\Gamma_h$. 
 \\
\noindent For $S \in {\mathcal S}_h$ we denote by ${\mathcal L}_S$ the set of $k-1$ Lagrangian nodes of $S$ different from vertexes that belong to $\Gamma_h$. We also denote by ${\bf n}_h$ the unit outer normal vector along $\Gamma_h$. Since $div \; {\bf b} \equiv 0$ by assumption, integration by parts easily yields 
$\int_{\Omega_h} ({\bf b} \cdot {\bf grad} \; w) w = \displaystyle \oint_{\Gamma_h} \frac{{\bf b} \cdot {\bf n}_h}{2} w^2 $.  
\begin{equation}
\label{ahwv}
\begin{array}{l} 
a_h(w,v) = \displaystyle \sum_{T \in {\mathcal T}_h} \int_T \nu |{\bf grad} \; w |^2 \\
- \displaystyle \sum_{S \in {\mathcal S}_h} \left\{ \int_S \left[ 
\nu {\bf grad} \; w \cdot {\bf grad} \; r_S(w) + ({\bf b} \cdot {\bf grad} \; w ) r_S(w) \right] - \displaystyle \int_{e_S} \frac{{\bf b} \cdot {\bf n}_h}{2} w^2 \right\},
\end{array} 
\end{equation}
\noindent where $e_S$ is the edge of $S$ contained in $\Gamma_h$ and $r_S(w) = \displaystyle \sum_{M \in {\mathcal L}_S} w(M) \varphi_M$, $\varphi_M$ being the canonical basis function of the space ${\mathcal P}_k(S)$ associated with Lagrangian node $M$. \\
\noindent Now from standard results it holds for two mesh independent constants $C_{\varphi,0}$ and $C_{\varphi,1}$:
\begin{equation}
\label{phiestim}
\left\{
\begin{array}{l}
\parallel \varphi_M \parallel_{0,S} \leq C_{\varphi,0} h_S^2, \\
\\
\parallel {\bf grad} \; \varphi_M \parallel_{0,S} \leq C_{\varphi,1}.
\end{array}
\right.
\end{equation} 
where $\parallel \cdot \parallel_{0,S}$ denotes the norm of $L^2(S)$. \\
On the other hand, since $w(P)=0$, where $P$ is the point of $\Gamma$ corresponding to $M \in \Gamma_h$ in accordance with the definition of $W_h$, a simple Taylor expansion about $P$ allows us to conclude that $|w(M)| \leq length(\overline{PM}) \parallel {\bf grad} \; w \parallel_{0,\infty,S}$. Hence, 
for a suitable constant $C_{\Gamma}$ independent of $S$ we have, $|w(Q)| \leq C_{\Gamma} h_S^2 \parallel {\bf grad} \; w \parallel_{0,\infty,S} \forall Q \in 
e_S$, where $\parallel \cdot \parallel_{0,\infty,S}$ is the standard norm of $L^{\infty}(S)$.
Notice that $w$ vanishes identically along a polynomial curve interpolating the $k+1$ points of the set consisting of the $k-1$ points $P$ on $\Gamma$ plus the two vertexes of $S$ belonging to $\Gamma$. Thus a similar argument allows us to extend to all $Q \in e_S$ this estimate of $w(M)$, there is
\begin{equation}
\label{wqestim}
|w(Q)| \leq C_{\Gamma} h_S^2 \parallel {\bf grad} \; w \parallel_{0,\infty,S} \forall Q \in e_S, 
\end{equation}
by eventually adjusting the constant $C_{\Gamma}$.\\
Let $l_S$ denote the length of $e_S$. Using \eqref{wqestim} we can assert that 
\begin{equation}
\label{ointegestim}
\int_{e_S} {\bf b} \cdot {\bf n}_h w^2 \leq C_{\Gamma}^2 h_S^4 l_S \parallel {\bf b} \parallel_{0,\infty} \parallel {\bf grad} \; w \parallel_{0,\infty,S}^2.
\end{equation}
Moreover, from a classical inverse inequality, we may write for another mesh independent constant $C_{\infty}$:
\begin{equation}
\label{inverseineq}
\parallel {\bf grad} \; w \parallel_{0,\infty,S} \leq C_{\infty} h_S^{-1} \parallel 
{\bf grad} \; w \parallel_{0,S}.
\end{equation}
Hence noticing that $l_S \leq h_S$ and $card({\mathcal L}_S) = k-1$ $\forall S$, plugging \eqref{phiestim}, \eqref{ointegestim} and \eqref{inverseineq} into (\ref{ahwv}), we easily derive:
\begin{equation}
\label{ahwvbound} 
\begin{array}{l}
a_h(w,v) \geq \int_{\Omega_h} \nu |{\bf grad} \; w |^2 \\
- \displaystyle \left[ \left( \frac{C_{\Gamma}^2 C_{\infty}^2}{2} + C_{\infty} C_{\Gamma} C_{\varphi,0} \right) \parallel {\bf b} \parallel_{0,\infty} h^3 +  
C_{\infty} C_{\Gamma} C_{\varphi,1} \nu h \right]   
\displaystyle \sum_{S \in {\mathcal S}_h} (k-1) \parallel {\bf grad} \; w \parallel_{0,S}^2. 
\end{array}
\end{equation}         
From \eqref{ahwvbound} we readily obtain for two suitable mesh independent constants $C_0$ and $C_1$:
\begin{equation}
\label{ahwvbelow} 
a_h(w,v) \geq [ \nu(1 - C_1 h) - C_0 \parallel {\bf b} \parallel_{0,\infty} h^3] \parallel {\bf grad} \; w \parallel_{0,\Omega_h}^2
\end{equation}
Now using arguments in all similar to those employed above, we easily infer that,  
\begin{equation}
\label{normbound}
 \parallel {\bf grad} \; v \parallel_{0,h} \leq \parallel {\bf grad} \; w  \parallel_{0,h} + \parallel {\bf grad}(v-w) \parallel_{0,h} \leq (1+ C_1 h) 
\parallel {\bf grad} \; w \parallel_{0,h}.
\end{equation}
Combining (\ref{ahwvbelow}) and (\ref{normbound}), provided $h \leq \displaystyle \min[(4C_1)^ {-1},(4 C_0 \mbox{P\' e})^{-1/3}]$, where 
P\'e := $\parallel {\bf b} \parallel_{0,\infty}/\nu$ is the P\'eclet number, we establish (\ref{infsup}) with $\alpha = 2 \nu / 5$. \\

\noindent Since obviously $dim(V_h) = dim(W_h)$, the simple fact that (\ref{infsup}) holds implies that (\ref{Poissonh}) is uniquely solvable (cf. \cite{COAM}). \QED \\

\indent We will also need 
\begin{corollary}
\label{corol1}
Provided $h$ is sufficiently there exists a constant $\alpha^{'} > 0$ independent of $h$ such that,
\begin{equation}
\label{infsuprime}
\forall w \in W_h \neq 0, \displaystyle \sup_{v \in V_h \setminus \{ 0 \}} \frac{a_h(w,v)}{\parallel {\bf grad} \; w \parallel_{0} \parallel {\bf grad} \; v \parallel_{0}} 
\geq \alpha^{'}.
\end{equation}  
\end{corollary}

\prov The proof of (\ref{infsuprime}) is a simple variant of the one of (\ref{infsup}) thanks to the observation that
\[ 
C_U \parallel {\bf grad} w \parallel_{0,T^{'}} \geq \parallel {\bf grad} w \parallel_{0,T} \geq C_L \parallel {\bf grad} w \parallel_{0,T^{'}} \; \forall T \in {\mathcal S}_h \mbox{ and } \forall w \in W_h, \] 
for suitable mesh independent constants $C_U$ and $C_L$, since $w$ is a polynomial in $T^{'}$ (cf. \cite{Strang}). \QED 

\section{Error estimates}

In order to derive error estimates for problem (\ref{Poissonh}) we resort to the approximation theory of non coercive linear variational problems (cf. \cite{Babuska}, 
\cite{Brezzi} and \cite{COAM}). At this point it is important to recall that since $d \equiv 0$, the solution $u$ of (\ref{Poisson}) satisfies 
$a(u,v) = F(v)$ $\forall v \in H^1_0(\Omega)$, where 
\begin{equation}
\label{aF}
a(w,v) := \int_{\Omega} [\nu {\bf grad} \; w \cdot {\bf grad}\; v + ({\bf b} \cdot {\bf grad} \; w ) v ] \mbox{ and } F(v) : = \int_{\Omega} fv.
\end{equation} 
Hence, owing to the construction of $V_h$, if $\Omega$ is convex $u$ also fulfills $a_h(u,v) = F_h(v) \; \forall v \in V_h$. In case $\Omega$ is not convex, 
we could extend $u$ by zero in $\Omega_h \setminus \Omega$, to define $a_h(u,v)$. However in this case there will be a non zero residual 
$a_h(u,v)-F_h(v)$ for $v \in V_h$ whose order may erode the one the approximation method (\ref{Poissonh}) is supposed to attain. Nevertheless in this case such an effect can be neutralized by means of a trick to be explained later on. For the moment let us assume that $\Omega$ is convex.\\
\indent Let us denote by $\parallel \cdot \parallel_{r,D}$ (resp. $| \cdot |_{r,D}$) the standard norm (resp. semi-norm) of Sobolev space $H^{r}(D)$ for $r >0$ (cf. \cite{Adams}), $D$ being any bounded 
domain of $\Re^2$ with non zero measure. We have,
\begin{theorem}
\label{theorem1}
As long as $h$ is sufficiently small, if $\Omega$ is convex and the solution $u$ of (\ref{Poisson}) for $d \equiv 0$ belongs to $H^{k+1}(\Omega)$, the solution $u_h$ of (\ref{Poissonh}) satisfies for $k>1$ and a suitable constant ${\mathcal C}$ independent of $h$ and $u$:
\begin{equation}
\label{errestconvex}
\parallel {\bf grad}(u - u_h) \parallel_{0,h} \leq {\mathcal C} h^k | u |_{k+1,\Omega}.
\end{equation}
\end{theorem}

\prov First we note that $u$ belongs to $H^2(\Omega)$. Therefore it is possible to uniquely define $u(Q)$ at any point $Q \in \Omega$ (cf. \cite{Adams}), and hence 
a $W_h$-interpolate of $u$ that we denote by $I_h(u)$. More specifically $I_h(u)$ is defined in the following fashion. In every $T \in 
{\mathcal T}_h \setminus {\mathcal S}_h$, $I_h(u)$ is the standard ${\mathcal P}_k$-interpolate of $u$ at the Lagrangian nodes of $T$. If $T \in {\mathcal S}_h$, 
$I_h(u)$ is the ${\mathcal P}_k$-interpolate of $u$ in $T^{'}$ at the set of $m_k+2$ points consisting of the Lagrangian nodes of $T$ that do not lie in the 
interior of $e_T$, together with the $k-1$ points $P$ lying on $\Gamma$ associated with the Lagrangian nodes $M$ 
of $T$ lying in the interior of $e_T$, as described in the definition of $W_h$. \\
From standard results (see e.g. \cite{Ciarlet}) we know that 
\begin{equation}
\label{errest2}
\parallel {\bf grad}[u - I_h(u)]_{|T} \parallel_{0,T} \leq C_{\Omega} h^k | u |_{k+1,T} \; \forall T \in {\mathcal T}_h \setminus {\mathcal S}_h,
\end{equation}
where $C_{\Omega}$ is a constant independent of $h$ and $u$. 
Now if $T \in {\mathcal S}_h$ we consider the mapping ${\mathcal G}_T$ from $T^{'}$ onto a unit element $\hat{T}^{'}$ of a reference plane with coordinates 
$(\hat{x},\hat{y})$ given by ${\mathcal G}_T(x,y) = (x,y)/h_T$. 
Since $\Delta_T$ is a small perturbation of $T$, $T^{'}$ is star-shaped with respect to a ball contained in $T$.
It follows that we can extend the well-known results for the Lagrange interpolation with the set of Lagrangian nodes to the one constructed in accordance with the definition of $W_h$. More precisely we mean the set consisting of the $m_k+2$ transformations in $\hat{T}$ under ${\mathcal G}_h$ of Lagrangian nodes of $T$  
which do not lie in the interior of $e_T$, completed with the transformations under ${\mathcal G}_h$ of the $k-1$ points $P \in \Gamma \cap T^{'}$ associated with the Lagrangian nodes $M$ of $T$ lying in the interior of $e_T$ (see Figure 2). Let us denote by $\hat{u}$ and $\widehat{I_h(u)}$ the transformations under ${\mathcal G}_h$ in $\hat{T}^{'}$ of $u$ and $I_h(u)$ restricted to $T^{'}$, respectively. 
Notice that $\widehat{I_h(u)}$ is the $P_k$-interpolate $\hat{I}(\hat{u})$ of $\hat{u}$ in $\hat{T}^{'}$, both functions coinciding 
whenever $\hat{u}$ belongs to ${\mathcal P}_k(\hat{T}^{'})$ and hence to ${\mathcal P}_k(\hat{T})$. Thus, denoting by $\rho_T$ the radius of the circle inscribed in $T$, 
by the same arguments as in Theorem 4.4.4 of \cite{BrennerScott}, we immediately conclude that 
\begin{equation}
\label{errest1}
 \parallel {\bf grad}[u - I_h(u)]_{|T^{'}} \parallel_{0,T^{'}} \leq C^{'} \displaystyle \frac{h_T^{k+1}}{\rho_T} | u |_{k+1,T^{'}} \; \forall T \in {\mathcal S}_h, 
\end{equation} 
\noindent $C^{'}$ being a constant depending only on $k$ and the set of interpolation points lying on $\Gamma \cap T^{'}$. Actually these points vary with $T$, but the underlying dependence of $C^{'}$ on them reduces to a dependence on $\Gamma$ rather than on $T$ itself.
\\
Now recalling that the chunkiness parameter $\sigma = \max_{T \in {\mathcal T}_h} h_T/\rho_T$ (cf. \cite{BrennerScott}) is bounded for every ${\mathcal T}_h$ in the family of partitions in use, we set $C_{\mathcal T} := C^{'} \sigma$.\\

On the other hand from (\ref{infsup}) we infer that

\begin{equation}
\label{bound1} 
\parallel {\bf grad}[u_h - I_h(u)] \parallel_{0,h} \leq \displaystyle \alpha^{-1} \sup_{v \in V_h \setminus \{0\}} \frac{a_h(u_h-I_h(u),v)}{\parallel {\bf grad} \; 
v \parallel_{0,h}} \; \forall w \in W_h.
\end{equation} 
Let us add and subtract $u$ in the first argument of $a_h$ and resort to the Friedrichs-Poincar\'e inequality, according to which $\parallel v \parallel_{0,h} 
\leq C_P \parallel {\bf grad}\; v \parallel_{0,h}$, where $C_P$ is constant depending only on $\Omega$. In doing so we obtain after a straightforward calculation:
\begin{equation}
\label{bound2} 
\parallel {\bf grad}[u_h - I_h(u)] \parallel_{0,h} \leq \displaystyle \alpha^{-1} \left[ A \parallel {\bf grad}[u - I_h(u)] \parallel_{0,h} + 
\displaystyle \sup_{v \in V_h \setminus \{0\}} \frac{a_h(u_h-u,v)}{\parallel {\bf grad} \; v \parallel_{0,h}} \right], 
\end{equation}
where $A := \nu + C_P \parallel {\bf b} \parallel_{0,\infty}$.
Noting that $a_h(u_h,v)=F_h(v)$ we come up with:
\begin{equation}
\label{bound3} 
\parallel {\bf grad}[u_h - I_h(u)] \parallel_{0,h} \leq \displaystyle \frac{1}{\alpha} \left\{ A \parallel {\bf grad}[u - I_h(u)] \parallel_{0,h} + 
\displaystyle \sup_{v \in V_h \setminus \{0\}} \frac{|a_h(u,v)-F_h(v)|}{\parallel {\bf grad}\; v \parallel_{0,h}} \right\}.
\end{equation}
Since $\Omega_h \subset \Omega$ if $\Omega$ is convex, we observe that $a_h(u,v) = \displaystyle \oint_{\Gamma_h} \nu v \displaystyle \frac{\partial u}{\partial n_h} +  \int_{\Omega_h} v (-\nu \Delta u + {\bf b} \cdot {\bf grad} \; u)$, where $\displaystyle \frac{\partial u}{\partial n_h}$ is the outer normal derivative of $u$ on $\Gamma_h$. From equation \eqref{Poisson} and since $v \equiv 0$ on $\Gamma_h$, it trivially follows that,
\begin{equation}
\label{bound4} 
\parallel {\bf grad}(u_h - u) \parallel_{0,h} \leq \displaystyle \left(1+ \frac{A}{\alpha} \right) \parallel {\bf grad}[u - I_h(u)] \parallel_{0,h}.
\end{equation}
Finally combining (\ref{bound4}), (\ref{errest2}) and (\ref{errest1}), we establish (\ref{errestconvex}) with 
${\mathcal C}:= [1+A/\alpha]C_{\mathcal T}$. \QED \\
\begin{corollary}
\label{convconvexprime}
As long as $h$ is sufficiently small, if $\Omega$ is convex and the solution $u$ of (\ref{Poisson}) for $d \equiv 0$ belongs to $H^{k+1}(\Omega)$, the solution $u_h$ of (\ref{Poissonh}) satisfies for $k>1$ and a suitable constant ${\mathcal C}^{'}$ independent of $h$ and $u$:
\begin{equation}
\label{errestconvexprime}
\parallel {\bf grad}(u - u_h) \parallel_{0} \leq {\mathcal C}^{'} h^k | u |_{k+1,\Omega}.
\end{equation}
\end{corollary}

\prov First we recall that the solution $u_h \in W_h$ of (\ref{Poissonh}) is also the solution of $a(u,v)=a(u_h,v)=F(v)$ $\forall v \in V_h$. Then owing to the fact that $a(u,v) \leq \parallel {\bf grad} \; u \parallel_0 \parallel {\bf grad} \; v \parallel_0$ $\forall (u,v) \in (H^1\Omega) + W_h) \times V_h$ and to (\ref{infsuprime}) we can write (cf. \cite{COAM}):
\[  \parallel {\bf grad}(u - u_h) \parallel_0 \leq \displaystyle \frac{A}{\alpha^{'}} \parallel {\bf grad}[u - I_h(u)] \parallel_0. \]
Then using (\ref{errest2}) and (\ref{errest1}) the result follows. \QED \\

\indent $O(h^{k+1})$-error estimates in the $L^2$-norm can be established in connection with 
Theorem \ref{theorem1}, if we require a little more regularity from $u$, according to,
\begin{theorem}
\label{theorem1bis}
As long as $h$ is sufficiently small, if $\Omega$ is convex and the solution $u$ of (\ref{Poisson}) for $g \equiv 0$ belongs to $H^{k+1+r}(\Omega)$ with $r=1/2+ \epsilon$ for $\epsilon > 0$ arbitrarily small, the solution $u_h$ of (\ref{Poissonh}) satisfies for $k>1$ and a suitable constant ${\mathcal C}_0$ independent of $h$ and $u$:
\begin{equation}
\label{L2estconvex}
\parallel u - u_h \parallel_{0,h} \leq {\mathcal C}_0 h^{k+1} \parallel u \parallel_{k+1+r,\Omega}.
\end{equation}
\end{theorem}

\prov Recalling that every function in $W_h$ is defined in $\Omega \setminus \Omega_h$, let $\bar{u}_h$ be the function given by $\bar{u}_h = u_h-u$ in $\Omega$. Let also $v \in H^1_0(\Omega)$ be the solution of 
\begin{equation}
\label{adjoint}
- \nu \Delta v - {\bf b} \cdot {\bf grad}\; v = \bar{u}_h \; \in \Omega.
\end{equation}
Since $\Omega$ is smooth and $\bar{u}_h \in L^2(\Omega)$ we know that $v \in H^2(\Omega)$ and moreover there 
exists a constant $C_{\Omega}$ depending only on $\nu$, ${\bf b}$ and $\Omega$ such that,
\begin{equation}
\label{adjoint1}
\parallel v \parallel_{2,\Omega} \leq C_{\Omega} \parallel \bar{u}_h \parallel_{0,\Omega}. 
\end{equation}
Therefore
\begin{equation}
\label{L2est1} 
\parallel \bar{u}_h \parallel_{0,\Omega} \leq C_{\Omega} \displaystyle \frac{\int_{\Omega} \bar{u}_h (- \nu \Delta v - {\bf b} \cdot {\bf grad}\; v)}{\parallel v \parallel_{2,\Omega}}.
\end{equation}
Using integration by parts we easily obtain,
\begin{equation}
\label{L2est2} 
\parallel \bar{u}_h \parallel_{0,\Omega} \leq C_{\Omega} \displaystyle \frac{a(\bar{u}_h, v)
+ b_{1h}(\bar{u}_h,v)}{\parallel v \parallel_{2,\Omega}} 
\end{equation}
where
\begin{equation}
\label{b1h}
b_{1h}(w,v):= -\nu \displaystyle \int_{\Gamma} w  \frac{\partial v}{\partial n} \mbox{ for } w \in H^1(\Omega) 
\mbox{ and } v \in H^1_0(\Omega).
\end{equation}
Let $\Pi_h(v)$ be the continuous piecewise linear interpolate of $v$ in $\Omega$ at the vertices of the mesh. 
Setting $v_h=\Pi_h(v)$ in $\Omega_h$ and $v_h=0$ in $\Omega \setminus \Omega_h$ we have $v_h \in V_h$. Therefore it holds $a(u,v_h)=a_h(u,v_h)=F(v_h)=F_h(v_h)=a_h(u_h,v_h)$. On the other hand 
$a(\bar{u}_h,v)=a_h(\bar{u}_h,v) + a_{\Delta_h}(\bar{u}_h,v)$ where 
\begin{equation}
\label{aDeltah}
a_{\Delta_h}(w,z):=\int_{\Delta_h} [\nu {\bf grad}\; w \cdot {\bf grad}\;z + {\bf b} \cdot {\bf grad}\; w \; z] \mbox{ for } w,z \in H^1(\Omega) \mbox{ with } \Delta_h = \Omega \setminus \Omega_h.
\end{equation}
Now we observe that $a_{\Delta_h}(\bar{u}_h,v)=a_{\Delta_h}(\bar{u}_h,v-\Pi_h(v))
+a_{\Delta_h}(\bar{u}_h,\Pi_h(v))$. Thus applying First Green's identity in $\Delta_T$ for $T \in {\mathcal S}_h$we come up with,
$a_{\Delta_h}(\bar{u}_h,\Pi_h(v)) = b_{2h}(\bar{u}_h,\Pi_h(v))+b_{3h}(\bar{u}_h,\Pi_h(v))$, where
\begin{equation}
\label{b2h}
b_{2h}(w,z):=\displaystyle \sum_{T \in {\mathcal S}_h} \int_{\Delta_T} [-\nu \Delta w+{\bf b} \cdot {\bf grad} \; w] z \mbox{ for } w \in W_h + H^2(\Omega) \mbox{ and } z \in H^{1}(\Omega),  
\end{equation}
and setting $\partial T = T^{'} \cap \Gamma$ for $T \in {\mathcal S}_h$,  
\begin{equation}
\label{b3h}
b_{3h}(w,z) := \nu \displaystyle \sum_{T \in {\mathcal S}_h} \int_{\partial T} \frac{\partial w}{\partial n} z \mbox{ for } w \in W_h + H^2(\Omega) \mbox{ and } z \in H^1(\Omega).
\end{equation}
Further setting
\begin{equation}
\label{b4h} 
b_{4h}(w,z) := a_{\Delta_h}(w,z) \mbox{ for } w,z \in H^1(\Omega),  
\end{equation}  
it follows that,
\begin{equation}
\label{L2est3} 
\begin{array}{l}
\parallel \bar{u}_h \parallel_{0,\Omega} \leq C_{\Omega} \displaystyle \frac{a_h(\bar{u}_h,e_h(v))
+ b_{1h}(\bar{u}_h,v)+b_{2h}(\bar{u}_h,\Pi_h(v))+b_{3h}(\bar{u}_h,\Pi_h(v))+
b_{4h}(\bar{u}_h,e_h(v))}{\parallel v \parallel_{2,\Omega}}, \\
\mbox{with } e_h(v) = v -\Pi_h(v). 
\end{array}
\end{equation} 
Now from classical results, for a mesh-independent constant $C_{V}$ it holds
\begin{equation}
\label{interP1}
\parallel {\bf grad} \; e_h(v) \parallel_{0,h} \leq \parallel {\bf grad} \; e_h(v) \parallel_{0} \leq C_{V} h 
| v |_{2,\Omega}.
\end{equation} 
Therefore, combining \eqref{bound2}, \eqref{errestconvex}, \eqref{interP1} and \eqref{L2est3}, and setting  
$\tilde{\mathcal C}_0=C_{\Omega} A C_{V} {\mathcal C}$, it holds,
\begin{equation}
\label{L2est4} 
\parallel \bar{u}_h \parallel_{0,\Omega} \leq \tilde{\mathcal C}_0 h^{k+1} | u |_{k+1,\Omega} + C_{\Omega} \displaystyle \frac{b_{1h}(\bar{u}_h,v)+b_{2h}(\bar{u}_h,\Pi_h(v))+b_{3h}(\bar{u}_h,\Pi_h(v))+
b_{4h}(\bar{u}_h,e_h(v))}{\parallel v \parallel_{2,\Omega}}.
\end{equation} 

Let us estimate $b_{ih}$ for $i=1,2,3,4$.\\
As for $b_{1h}$ we first note that according to the Trace Theorem there exists a constant $C_t$ depending only on $\Omega$ such that 
\begin{equation}
\label{L2estimate}
b_{1h}(\bar{u}_h,v) \leq C_t \parallel \bar{u}_h \parallel_{0,\Gamma} \parallel v \parallel_{2,\Omega}.
\end{equation} 
Now for every $T \in {\mathcal S}_h$ we take a local orthogonal frame $(O;x,y)$ whose origin $O$ is a vertex of $T$ in $\Gamma$, $x$ is the abscissa along the edge $e_T$ and $y$ increases from 
$e_T$ towards $\Gamma$. Let $\partial T:=\Gamma \cap T^{'}$ and $s$ be the curvilinear abscissa along $\partial T$ with origin at $O$. Notice that owing to our assumptions $s$ can be uniquely expressed in terms of $x$ and conversely, for $x \in [0,l_T]$, where $l_T$ is the length of $e_T$. Let then $\breve{u}_h$ be the function of $x$ defined by $\breve{u}_h(x) = \bar{u}_h[s(x)]$. Since $\breve{u}_h$ vanishes at $k+1$ different points in $[0,l_T]$, from standard results for one-dimensional interpolation (cf \cite{Quarteroni}), there exists a mesh-independent constant $C_e$ such that,
\begin{equation}
\label{auxiliary0}
\displaystyle \left[\int_0^{l_T} |\breve{u}_h(x)|^2 dx \right]^{1/2} \leq C_e h_T^{k+1} \displaystyle \left[\int_0^{l_T} \left|\frac{d^{k+1}\breve{u}_h}{dx^{k+1}}(x) \right|^2 dx \right]^{1/2} 
\end{equation}
On the other hand defining the function $y(x)$ to be the $y$-abscissa of the points in 
$\partial T$, we observe that there exist mesh-independent constants $c_{j,\Gamma}$ such that,
\begin{equation}
\label{auxiliary1}
\displaystyle \max_{x \in [0,l_T]}  |y^{(j)}(x)| \leq c_{j,\Gamma} h_T^{2-j}, \; j=1,2,\ldots,k+1 \; \forall T \in {\mathcal S}_h.
\end{equation}
Thus taking into account that the derivatives of $u_h$ of order greater than $k$ vanish in $T$, straightforward calculations using the chain rule yield for suitable mesh-independent constants $c_{j}$, $j=0,1,\ldots,k$: \\
\begin{equation}
\label{auxiliary2}
\displaystyle \frac{d^{k+1}\breve{u}_h}{dx^{k+1}} \leq  c_{0} | D^{k+1}(u) | + \displaystyle 
\sum_{j=1}^{k} c_{j} h_T^{1-j} | D^{k+1-j}(\bar{u}_h) |,  
\end{equation}
where $D^j w$ is the $j$-th order tensor, whose components are the $j$-th order partial derivatives 
of a function $w$ in $\Omega$. \\ 
All the partial derivatives appearing in \eqref{auxiliary2} are to be understood at a (variable) point in $\partial T$. \\
Now since $ds = \sqrt{1+(y^{'})^2}dx$, there exists another mesh-independent constant $C_q$ such that 
\begin{equation}
\label{auxiliary3}
\parallel \bar{u}_h \parallel_{0,\Gamma} \leq 
C_q \displaystyle \left[ \sum_{T \in {\bf S}_h} \int_0^{l_T} |\breve{u}_h(x)|^2 dx \right]^{1/2}.
\end{equation}
Combining \eqref{auxiliary0}, \eqref{auxiliary1}, \eqref{auxiliary2} and \eqref{auxiliary3}, after straightforward calculations we come up with a mesh-independent constant $\tilde{C}$ such that,
\begin{equation}
\label{auxiliary4}
\parallel u_h \parallel_{0,\Gamma} 
\leq \tilde{C}  \displaystyle \left\{ \sum_{T \in {\bf S}_h} h_T^{2(k+1)} 
\displaystyle \int_{\partial T} \left[ | D^{k+1}(u) |^2 + \displaystyle 
\sum_{j=1}^k h_T^{2(1-j)}|D^{k+1-j}(\bar{u}_h)|^2  \right] \right\}^{1/2}.\\
\end{equation}
From the Trace Theorem \cite{Adams} we know that there exists a constant $C_{r}(\Omega)$ such that,
\begin{equation}
\label{auxiliary5}
\left[\displaystyle \sum_{T \in {\bf S}_h} 
\int_{\partial T}  |D^{k+1}(u)|^2 \right]^{1/2} \leq 
C_{r}(\Omega) \parallel u \parallel_{k+1+r,\Omega} 
\end{equation}
On the other hand, using the curved triangle $T^{'}$ associated with
$T$, by standard calculations we can write for a suitable mesh-independent constant $C_{k,1}$:
\begin{equation}
\label{auxiliary6}
\begin{array}{l}
\displaystyle \int_{\partial T}  
\displaystyle \left[ \sum_{j=1}^k h_T^{2(1-j)}|D^{k+1-j}(\bar{u}_h)|^2  \right]\\
\leq C_{k,1} h_T \displaystyle \sum_{j=1}^k h_T^{2(1-j)} 
\left[\parallel D^{k+1-j}(u_h-I_h(u)) \parallel_{0,\infty,T^{'}}^2 + 
\parallel D^{k+1-j}(I_h(u)-u) \parallel_{0,\infty,T^{'}}^2 \right].
\end{array}
\end{equation} 
Since $[D^{k+1-j}(u_h-I_h(u))]_{|T^{'}} \in {\mathcal P}_{j-1}$ and $area(\Delta_T)/area(T) = O(h_T)$ we 
have for suitable constants $C_{k}^j$ independent of $T$, $u_h-I_h(u)$:
\begin{equation}
\label{skin}
\parallel D^{k+1-j}(u_h-I_h(u)) \parallel_{0,\infty,T^{'}} \leq C_{k}^j \parallel D^{k+1-j}(u_h-I_h(u)) \parallel_{0,\infty,T}, \; j=1,2,\ldots,k.
\end{equation}
Noting that $H^{j+1+r}(\Omega)$ is embedded in $W^{j,\infty}(\Omega)$ for $j=0,1,\ldots,k$, the embedding being continuous (cf. \cite{Adams}), applying classical inverse inequalities in $T$, together with 
well-know estimates for the interpolation error, for another mesh-independent constant $C_{k,2}$ it holds: 
\begin{equation}
\label{auxiliary7}
\displaystyle \int_{\partial T}  
\left[ \sum_{j=1}^k h_T^{2(1-j)}|D^{k+1-j}(\bar{u}_h)|^2  \right] \\   
\leq \displaystyle  C_{k,2} \left\{ h_T^{-2k+1} \parallel {\bf grad}(u_h-I_h(u))  
\parallel_{0,T}^2 + h_T \parallel u \parallel_{k,\infty,\Omega}^2 \right\}, 
\end{equation}
where $\parallel \cdot \parallel_{l,\infty,D}$ denotes the standard norm of $W^{l,\infty}(D)$ for an 
integer $l>0$.\\ 
Now \eqref{errestconvex} together with the estimate $\parallel {\bf grad}(u-I_h(u)) \parallel_{0,T} \leq C_{k,3} h_T^{k} | u |_{k+1,T}$ for a suitable $C_{k,3}$ easily yield:
\begin{equation}
\label{auxiliary}
\displaystyle \sum_{T \in {\mathcal S}_h} h_T^{-2k+1} \parallel {\bf grad}(u_h-I_h(u)) \parallel_{0,T}^2 
\leq 2(C_{k,3}^2 + {\mathcal C}^2) h | u |_{k+1,\Omega}^2.
\end{equation} 
On the other hand we observe that by the Sobolev Embedding Theorem there exists a constant $C_E$ such that 
\begin{equation}
\label{embedding}
\parallel u \parallel_{k,\infty,\Omega} \leq C_E \parallel u \parallel_{k+1+r,\Omega}.
\end{equation}
Therefore for a certain constant ${\mathcal C}_{\Omega}$ we have  
\begin{equation}
\label{auxiliary8}
\displaystyle \sum_{T \in {\bf S}_h} h_T \parallel u \parallel_{k,\infty,\Omega}^2  \leq C_E^2   
\displaystyle \left[ \sum_{T \in {\bf S}_h} h_T \right] \parallel u \parallel_{k+1+r,\Omega}^2 
\leq {\mathcal C}_{\Omega} \parallel u \parallel_{k+1+r,\Omega}^2 
\end{equation}
Taking into account \eqref{auxiliary4}, \eqref{auxiliary5}, \eqref{auxiliary7}, \eqref{auxiliary} and \eqref{auxiliary8}, we easily obtain,
\begin{equation}
\label{auxiliary9}
\parallel \bar{u}_h \parallel_{0,\Gamma} 
\leq \bar{C}_1 h^{k+1} \displaystyle \left[h^{1/2} | u |_{k+1,\Omega} + \parallel u \parallel_{k+1+r,\Omega} \right],
\end{equation} 
where $\bar{C}_1$ is a mesh-independent constant.\\
It follows from \eqref{L2estimate} and \eqref{auxiliary9} that for $C_{b1} = 2 \bar{C}_1 C_t $ it holds:
\begin{equation}
\label{estimateb1}
b_{1h}(\bar{u}_h,v) \leq C_{b1} h^{k+1}[ h^{1/2} | u |_{k+1,\Omega }\parallel u \parallel_{k+1+r,\Omega}] \parallel v \parallel_{2,\Omega}.
\end{equation} 

Now we turn our attention to $b_{2h}$.\\
First observing that ${\bf grad}\;\Pi_h(v)$ is constant in $T^{'}$ for $T \in {\mathcal T}_h$ and $\Pi_h(v) =0$ on $\Gamma_h$, by Rolle's Theorem
\begin{equation}
\label{estim1b2}
|\Pi_h(v)(P)|\leq C_{\Gamma} h_T^2 \parallel {\bf grad} \; \Pi_h(v) \parallel_{0,\infty,T} \; \forall P \in \partial T \mbox{ and } \forall T \in {\mathcal S}_h.
\end{equation}
Noticing that $area(\Delta_T) \leq C_{\Gamma} h_T^3$, using \eqref{estim1b2} we have,
\begin{equation}
\label{estim2b2}
b_{2h}(\bar{u}_h,\Pi_h{v}) \leq C_{\Gamma}^2 \displaystyle \sum_{T \in {\mathcal S}_h}  h_T^5 
\parallel -\nu \Delta \bar{u}_h + {\bf b} \cdot {\bf grad} \; \bar{u}_h \parallel_{0,\infty,T^{'}} 
\parallel {\bf grad} \; \Pi_h(v) \parallel_{0,\infty,T}.
\end{equation}
Using the classical inverse inequality $\parallel {\bf grad} \; \Pi_h(v) \parallel_{0,\infty,T} \leq C_I h_T^{-1} \parallel {\bf grad} \; \Pi_h(v) \parallel_{0,T}$ with $C_I$ independent of $T$, we further obtain:
\begin{equation}
\label{estim3b2}
b_{2h}(\bar{u}_h,\Pi_h{v}) \leq C_{\Gamma}^2 C_I \displaystyle \sum_{T \in {\mathcal S}_h}  h_T^4 
\parallel -\nu \Delta \bar{u}_h + {\bf b} \cdot {\bf grad} \; \bar{u}_h \parallel_{0,\infty,T^{'}} 
\parallel {\bf grad} \; \Pi_h(v) \parallel_{0,T}.
\end{equation}
Next using the triangle inequality we rewrite \eqref{estim3b2} as, 
\begin{equation}
\label{estim4b2} 
\begin{array}{l}
b_{2h}(\bar{u}_h,\Pi_h{v}) \leq C_{\Gamma}^2 C_I \displaystyle \sum_{T \in {\mathcal S}_h}  h_T^4 
[ \parallel -\nu \Delta(u_h-I_h(u)) + {\bf b} \cdot {\bf grad}(u_h-I_h(u)) \parallel_{0,\infty,T^{'}} \\
+ \parallel -\nu \Delta(I_h(u)-u) + {\bf b} \cdot {\bf grad}(I_h(u)-u) \parallel_{0,\infty,T^{'}}]
\parallel {\bf grad} \; \Pi_h(v) \parallel_{0,T}.
\end{array}
\end{equation}
From the inverse inequality $\parallel \Delta(u_h-I_h(u)) \parallel_{0,\infty,T^{'}} 
\leq C^{'}_I  h_T^{-1} \parallel {\bf grad}(u_h-I_h(u)) \parallel_{0,\infty,T^{'}}$ for another constant $C^{'}_I$ independent of $T$, and again the above one, we have 
\begin{equation}
\label{estim5b2} 
\parallel -\nu \Delta(u_h-I_h(u)) + {\bf b} \cdot {\bf grad}(u_h-I_h(u)) \parallel_{0,\infty,T^{'}} 
\leq C_2^{'} h_T^{-2} \parallel {\bf grad}(u_h-I_h(u)) \parallel_{0,T^{'}} 
\end{equation}
where $C_2^{'}$ is a mesh-independent constant.\\
Plugging \eqref{estim5b2} into \eqref{estim4b2} and further using the triangle inequality, we easily obtain,
\begin{equation}
\label{estim6b2} 
\begin{array}{l}
b_{2h}(\bar{u}_h,\Pi_h{v}) \leq \tilde{C}_2 \displaystyle \sum_{T \in {\mathcal S}_h}\{ h_T^2 
[\parallel {\bf grad}(u_h-u) \parallel_{0,T^{'}} + \parallel {\bf grad}(u-I_h(u)) \parallel_{0,T^{'}}] \\
+ h_T^4 \sqrt{2} \parallel u-I_h(u) \parallel_{2,\infty,T^{'}} \}\parallel {\bf grad} \; \Pi_h(v) \parallel_{0,T},
\end{array}
\end{equation}
for a suitable mesh-independent constant $\tilde{C}_2$.\\
Using the Cauchy-Schwarz inequality and taking into account that $\displaystyle \sum_{T \in {\mathcal S}_h} h_T \leq C^{'}(\Gamma)$ where $ C^{'}(\Gamma)$ is a mesh-independent constant, from  
\eqref{estim6b2} we derive for another mesh-independent constant $\bar{C}_2$, 
\begin{equation}
\label{estim7b2} 
\begin{array}{l}
b_{2h}(\bar{u}_h,\Pi_h{v}) \leq \bar{C}_2 \{ h^2 [\parallel {\bf grad}(u_h-u) \parallel_{0,h} + 
\parallel {\bf grad}(u-I_h(u)) \parallel_{0,h} ] \\
+ h^{7/2}  \parallel u-I_h(u) \parallel_{2,\infty,\Omega} \}\parallel {\bf grad} \; \Pi_h(v) \parallel_{0,h}.
\end{array}
\end{equation}
From standard interpolation results and \eqref{embedding} we can assert that for three mesh-independent constants $C_{1,2}$, $C_{2,2}$ and $C_{3,2}$ it holds, 
\begin{equation}
\label{estim8b2}
\left\{
\begin{array}{l}
\parallel u-I_h(u) \parallel_{2,\infty,\Omega} \leq C_{1,2} h^{k-2} | u |_{k,\infty,\Omega} \leq 
C_{1,2} C_E h^{k-2} \parallel u \parallel_{k+1+r,\Omega}; \\
\\
\parallel {\bf grad}(u-I_h(u)) \parallel_{0,h} \leq C_{2,2} h^k | u |_{k+1,\Omega}; \\
\\
\parallel {\bf grad} \; \Pi_h(v) \parallel_{0,h} \leq C_{3,2} \parallel v \parallel_{2,\Omega}.
\end{array}
\right.
\end{equation}
Plugging \eqref{estim8b2} into \eqref{estim7b2}, and recalling \eqref{errestconvex} we finally obtain,
\begin{equation}
\label{estimateb2}
b_{2h}(\bar{u}_h,\Pi_h{v}) \leq C_{b2} h^{k+1}[ h | u |_{k+1,\Omega} +h^{1/2} \parallel u \parallel_{k+1+r,\Omega} ] \parallel v \parallel_{2,\Omega},
\end{equation} 
where $C_{b2}$ is a mesh-independent constant.\\

Next we estimate $b_{3h}$.\\
Recalling \eqref{b3h} and the fact $\parallel {\bf grad}\;\Pi_h(v) \parallel_{0,\infty,T^{'}} = \parallel {\bf grad}\;\Pi_h(v) \parallel_{0,\infty,T}$, we first define the function $\omega_T:=|{\bf grad}\;\bar{u}_{h_{|T}}|$ 
for every $T \in {\mathcal S}_h$. Then we have:
\begin{equation}
\label{estim1b3}
b_{3h}(\bar{u}_h,\Pi_h(v)) \leq \nu \displaystyle \sum_{T \in {\mathcal S}_h} 
\int_{\partial T} \omega_T \Pi_h(v) \leq \nu C_{\Gamma} h_T^2 \parallel {\bf grad} \; \Pi_h(v) \parallel_{0,\infty,T} \int_{\tilde{\partial} T} \omega_T.
\end{equation} 
Let us denote the standard master triangle by $\hat{T}$ and the transformation of $\partial T$ 
under the affine mapping ${\mathcal F}_T$ from $T$ onto $\hat{T}$ by $\hat{\partial}\hat{T}$. Clearly enough  there exists a constant $\hat{C}$ independent of $T$ such that,
\begin{equation}
\label{estim2b3}
b_{3h}(\bar{u}_h,\Pi_h(v)) \leq \nu C_{\Gamma} \hat{C} \displaystyle \sum_{T \in {\mathcal S}_h} h_T^3 
\parallel {\bf grad} \; \Pi_h(v) \parallel_{0,\infty,T}  \int_{\hat{\partial} \hat{T}} \hat{\omega}, 
\end{equation}  
where $\hat{\omega}$ is the transformation of $\omega_T$ under the mapping 
${\mathcal F}_T$. We denote by $\hat{T}^{'}$ the transformation of $T^{'}$ under ${\mathcal F}_T$.\\
Next we apply the Trace Theorem to $\hat{T}^{'}$. Thanks to the fact that $\Gamma$ is smooth and $h$ is sufficiently small, there exists a constant $\hat{C}_t$ independent of $T$ such that,
\begin{equation}
\label{estim3b3} 
\displaystyle \int_{\hat{\partial}\hat{T}} \hat{\omega} \leq \hat{C}_t \displaystyle \left\{ \int_{\hat{T}^{'}} [ \hat{\omega}^2 + |\widehat{\bf grad}\; \hat{\omega}|^2 ]  \right\}^{1/2},   
\end{equation}  
where $\widehat{\bf grad}$ is the gradient operator for functions defined in $\hat{T}^{'}$.\\
Moving back to $T^{'}$ associated with $T \in {\mathcal S}_h$ and using an inverse inequality, from \eqref{estim2b3} and \eqref{estim3b3} we conclude that for a suitable mesh-independent constant $\breve{C}_3$ it holds,
\begin{equation}
\label{estim4b3}
b_{3h}(\bar{u}_h,\Pi_h(v)) \leq \breve{C}_{3} \displaystyle \sum_{T \in {\mathcal S}_h} 
h_T \parallel {\bf grad} \; \Pi_h(v) \parallel_{0,T} \left\{ \int_{T^{'}} [ \omega_T^2 + h_T^2|{\bf grad}\;\omega_T|^2 ] \right\}^{1/2}. 
\end{equation}
By the Cauchy-Schwarz inequality this further yields,  
\begin{equation}
\label{estim5b3}
b_{3h}(\bar{u}_h,\Pi_h(v)) \leq \breve{C}_{3} h \parallel {\bf grad} \; \Pi_h(v) \parallel_{0,h} 
[ \parallel {\bf grad}\;\bar{u}_{h} \parallel_{0}^2  + h^2 \parallel H(\bar{u}_{h})  
\parallel_{0}^2 ]^{1/2}. 
\end{equation}
Now using the triangle inequality and an inverse inequality, and then combining \eqref{errestconvex} with the second equation of \eqref{estim8b2}, we can estimate $\parallel H(\bar{u}_h) \parallel_{0,h}$ 
in the same way as $\parallel \Delta \bar{u}_h \parallel_{0,h}$ starting from \eqref{estim4b2}. In this way we can easily establish the existence of a mesh-independent constant $\bar{C}_3$ such that, 
\begin{equation}
\label{estim6b3}
[ \parallel {\bf grad}\;\bar{u}_{h} \parallel_{0}^2  + h^2 \parallel H(\bar{u}_{h})  
\parallel_{0}^2 ]^{1/2} \leq \bar{C}_3 h^{k} | u |_{k+1,\Omega}.
\end{equation}
Then plugging the third equation of \eqref{estim8b2} and \eqref{estim6b3} into \eqref{estim5b3}, 
yields for $C_{b3}=\bar{C}_3 \breve{C}_{3} C_{2,3}$:
\begin{equation}
\label{estimateb3} 
b_{3h}(\bar{u}_h,\Pi_h(v)) \leq C_{b3} h^{k+1} | u |_{k+1,\Omega} \parallel v \parallel_{2,\Omega}. 
\end{equation}

Finally we estimate $b_{4h}$.\\
Using a few arguments already exploited above, we can write:
\begin{equation}
\label{estim1b4}
b_{4h}(\bar{u}_h,v-\Pi_h(v)) \leq A C_{\Gamma}^{1/2} \displaystyle \sum_{T \in {\mathcal S}_h} h_T^{3/2} 
\parallel {\bf grad} \; \bar{u}_h \parallel_{0,\infty,T^{'}} \parallel {\bf grad}( v - \Pi_h(v)) 
\parallel_{0,T^{'}}, 
\end{equation}
and further,
\begin{equation}
\label{estim2b4}
\begin{array}{l}
b_{4h}(\bar{u}_h,v-\Pi_h(v)) \leq A C_{\Gamma}^{1/2} \displaystyle \sum_{T \in {\mathcal S}_h} h_T^{3/2} 
[C_I h_T^{-1} \parallel {\bf grad}(u_h-I_h(u)) \parallel_{0,T^{'}} \\
+ \parallel {\bf grad}(I_h(u)-u) \parallel_{0,\infty,T^{'}}]
\parallel {\bf grad}( v - \Pi_h(v)) \parallel_{0,T^{'}}, 
\end{array}
\end{equation}

\begin{equation}
\label{estim3b4}
\begin{array}{l}
b_{4h}(\bar{u}_h,v-\Pi_h(v)) \leq A C_{\Gamma}^{1/2} \displaystyle \sum_{T \in {\mathcal S}_h}  
\{ C_I h_T^{1/2} [\parallel {\bf grad}(u_h-u) \parallel_{0,T^{'}} + 
\parallel {\bf grad}(I_h(u)-u) \parallel_{0,T^{'}}] \\
+ h_T^{3/2} \parallel {\bf grad}(I_h(u)-u) \parallel_{0,\infty,T^{'}} \} 
\parallel {\bf grad}( v - \Pi_h(v)) \parallel_{0,T^{'}}. 
\end{array}
\end{equation}
Since $\parallel {\bf grad}(I_h(u)-u) \parallel_{0,\infty,T^{'}}$ can be bounded above by a constant 
independent of $T^{'}$ multiplied by $h^{k-1} | u |_{k,\infty,\Omega}$, from the Cauchy-Schwarz inequality and using \eqref{interP1} together with \eqref{errestconvex}, we infer the existence of a mesh-independent constant $\breve{C}_{4}$ such that,
\begin{equation}
\label{estim4b4}
b_{4h}(\bar{u}_h,v-\Pi_h(v)) \leq \breve{C}_{4} \displaystyle \left\{ h^{k+3/2} | u |_{k+1,\Omega} 
+ \displaystyle \left[ \sum_{T \in {\mathcal S}_h} h_T \right]^{1/2} h^{k+1} | u |_{k,\infty,\Omega} \right\}
| v |_{2,\Omega}. 
\end{equation}
Taking into account \eqref{embedding} this implies in turn that for another mesh-independent constant 
$C_{b4}$ it holds,
\begin{equation}
\label{estimateb4}
b_{4h}(\bar{u}_h,v-\Pi_h(v)) \leq C_{b4} h^{k+1}[ h^{1/2} | u |_{k+1,\Omega} 
+ \parallel u \parallel_{k+r,\Omega} ] | v |_{2,\Omega}. 
\end{equation}
 
Plugging \eqref{estimateb1}, \eqref{estimateb2}, \eqref{estimateb3} and \eqref{estimateb4} into 
\eqref{L2est4}, owing to the fact that $h <1$, we immediately obtain \eqref{L2estconvex} with 
${\mathcal C}_0=\tilde{\mathcal C}_0 + 2 C(\Omega) (C_{b1}+C_{b2}+C_{b3}+C_{b4})$. \QED \\

Now we address the case of a non convex $\Omega$. Let us consider a smooth domain $\tilde{\Omega}$ close to $\Omega$ which strictly contains $\Omega \cup \Omega_h$ for all $h$ sufficiently small. More precisely, denoting by $\tilde{\Gamma}$ the boundary of $\tilde{\Omega}$ we assume that $meas(\tilde{\Gamma})-meas(\Gamma) \leq \varepsilon$ for $\varepsilon$ sufficiently small. Henceforth we also consider that $f$ was extended to  
$\tilde{\Omega} \setminus \Omega$. We still denote the extended function by $f$, which is arbitrarily chosen in $\tilde{\Omega} \setminus \Omega$, except for the requirement that $f \in H^{k-1}(\tilde{\Omega})$. 

Then under the conditions specified therein the following theorem holds:

\begin{theorem}
\label{theorem2}
Assume that there exists a function $\tilde{u}$ defined in 
$\tilde{\Omega}$ having the properties:
\begin{itemize}
\item 
$- \nu \Delta \tilde{u} + {\bf b} \cdot {\bf grad} \; \tilde u = f$ in $\tilde{\Omega}$;
\item
$\tilde{u}_{|\Omega} = u$;
\item
$\tilde{u} = 0$ a.e. on $\Gamma$; 
\item
$\tilde{u} \in H^{k+1}(\tilde{\Omega})$.
\end{itemize}
Then as long as $h$ is sufficiently small it holds:   
\begin{equation}
\label{errestconcave}
\parallel {\bf grad}(u_h - u) \parallel_{0,\tilde{\Omega}_h} \leq \tilde{\mathcal C} h^k | \tilde{u} |_{k+1,\tilde{\Omega}},
\end{equation}
where $\tilde{\mathcal C}$ is a mesh-independent constant and $\parallel \cdot \parallel_{0,\tilde{\Omega}_h}$ denotes the standard $L^2$-semi-norm in the set 
$\tilde{\Omega}_h := \Omega_h \cap \Omega$. 
\end{theorem}

\prov Thanks to its properties $\tilde{u}$ can replace $u$ in the proof of Theorem \ref{theorem1}, to transform it into the proof of this theorem based on the same arguments. Then the observation that $\parallel {\bf grad}(u - u_h) \parallel_{0,\tilde{\Omega}_h} \;
\leq \; \parallel {\bf grad}(\tilde{u} - u_h) \parallel_{0,h}$ leads to (\ref{errestconcave}). \QED \\

It is noteworthy that the knowledge of a regular extension of the right hand side datum $f$ associated with a regular extension $\tilde{u}$ of $u$ is necessary to 
optimally solve problem (\ref{Poissonh}) in the general case. Of course, except for very particular situations such as the toy problems used to illustrate the performance of our method in the next section, in most cases such an extension of $f$ is not known. Even if we go the other around by prescribing a regular $f$ in 
$\tilde{\Omega}$, the existence of an associated $\tilde{u}$ satisfying the assumptions of Theorem \ref{theorem2} can also be questioned. However using some results available in the literature it is possible to identify cases where such an extension $\tilde{u}$ does exist. Let us consider for instance the Poisson equation (that is, $\nu=1$ and ${\bf b} \equiv {\bf 0}$) 
in a simply connected domain $\Omega$ of the $C^{\infty}$-class and a datum $f$ 
infinitely differentiable in $\bar{\Omega}$. Taking an extension of $f$ to the enlarged domain 
$\tilde{\Omega}$  also of the $C^{\infty}$-class, such that $f \in C^{\infty}(\tilde{\Omega}) \cap H^{k-1}(\tilde{\Omega})$, we first solve $-\Delta u_0 = f$ in 
$\tilde{\Omega}$ and $u_0 = 0$ on $\tilde{\Gamma}$. According to well-known results (cf. \cite{LionsMagenes}) $u_0 \in C^{\infty}(\tilde{\Omega})$ and hence the trace $g_0$ of $u_0$ on $\Gamma$ belongs to $C^{\infty}(\Gamma)$. Next we denote by $u_H$ the harmonic function in $\Omega$ such that $u_H = g_0$ on $\Gamma$. Let $r_0$ be the radius of the largest (open) ball $B$ contained in $\Omega$ and $O=(x_0,y_0)$ be its center. Assuming that the extension of $f$ is not too wild in $\tilde{\Omega}$ so that the Taylor series of $u_H(x,y_0)$ and 
$[\partial u_H/\partial y](x,y_0)$ centered at $O$ converge in the segment of the line $y=y_0$ centered at $O$ with length equal to 
$r_0 \sqrt{2} r_0 + 2\delta$ for a certain $\delta > 0$, according to \cite{Coffman} there exists a harmonic extension $u^{'}_H$ of $u_H$ to the ball $B_0^{'}$ 
centered at $O$ with radius $r_0 +\delta \sqrt{2}$. Clearly in this case, as long as $\delta$ is large enough for $B^{'}$ to contain $\tilde{\Omega}$, we can define 
$\tilde{u}:=u_0-u_H^{'}$ as a function in $H^{k+1}(\tilde{\Omega})$ that vanishes on $\Gamma$, and thus satisfies all the required properties.\\  
In the general case however, a convenient way to bypass the uncertain existence of an extension $\tilde{u}$ satisfying the assumptions of Theorem \ref{theorem2}, is to resort to numerical integration on the right hand side. Under certain conditions rather easily satisfied, this leads to the definition of an alternative approximate problem, in which only values of $f$ in $\Omega$ come into play. This trick is inspired by the celebrated work due to Ciarlet and Raviart on the isoparametric finite element method (cf. \cite{CiarletRaviart} and \cite{Ciarlet}). To be more specific, these authors employ the following argument, assuming that $h$ is small enough: if a numerical integration formula is used, which has no integration points different from vertexes on the edges of a triangle, then only values of $f$ in $\Omega$ will be needed to compute the corresponding approximation of $F_h(v)$. This means that the knowledge of $\tilde{u}$, and thus of the regular extension of $f$, will not be necessary for implementation purposes. Moreover, provided the accuracy of the numerical integration formula is compatible with method's order, the resulting modification of (\ref{Poissonh}) will be a method of order $k$ in the norm  $\parallel \cdot \parallel_{0,\tilde{\Omega}_h}$ of ${\bf grad}(u - u_h)$. \\
Nevertheless it is possible to get rid of the above argument based on numerical integration in the most important cases in practice, namely, the one of quadratic and cubic Lagrange finite elements. Let us see how this works. \\
First of all we consider that $f$ is extended by zero in $\Delta_{\Omega}:=\tilde{\Omega} \setminus \bar{\Omega}$, and resort to the extension $\tilde{u}$ of $u$ to the same set constructed in accordance to Stein et al. \cite{Stein}. This extension does not satisfy $\Delta \tilde{u}=0$ in $\Delta_{\Omega}$ but the function denoted in the same way such that $\tilde{u}_{|\Omega} = u$ does belong to $H^{k+1}(\tilde{\Omega})$. Since $k > 1$ this means in particular that the traces of the functions $u$ and $\tilde{u}$ coincide on $\Gamma$ and that $\partial u /\partial n = - \partial \tilde{u}/\partial \tilde{n} = 0$ a.e. on $\Gamma$ where the normal derivatives on the left and right hand side of this relation are outer normal derivatives with respect to $\Omega$ and $\Delta_{\Omega}$ respectively (the trace of the Laplacian of both functions also coincide  on $\Gamma$ but this is not relevant for our purposes). Based on this extension of $u$ to $\Omega_h$ for all such polygons of interest, we next prove the following results for the approximate problem (\ref{Poissonh}), without assuming that $\Omega$ is convex. Here $f$ represents the function identical in $\Omega$ to the right hand side datum of (\ref{Poisson}), that vanishes identically in $\Delta_{\Omega}$.

\begin{theorem} 
\label{P2}
Let $k=2$ and assume that $u \in H^3(\Omega)$. Provided $h$ is sufficiently small, there exists a mesh independent constant $C_2$ such that the unique solution $u_h$ to (\ref{Poissonh}) satisfies:
\begin{equation}
\label{estimateP2} 
\parallel {\bf grad}(u - u_h) \parallel_{0,\tilde{\Omega}_h} \leq C_2 [h^{2} | \tilde{u} |_{3,\tilde{\Omega}} + h^{5/2} \parallel \nu \Delta \tilde{u} 
- {\bf b} \cdot {\bf grad} \; \tilde{u} \parallel_{0,\tilde{\Omega}}]
\end{equation}
where $\tilde{u} \in H^3(\tilde{\Omega})$ is the regular extension of $u$ to $\tilde{\Omega}$ constructed in accordance to Stein et al. \cite{Stein}.
\end{theorem}

\prov
First we note that, 
\begin{equation}
\label{fourthbound}
\parallel {\bf grad}(u_h - w) \parallel_{0,h} \leq \displaystyle \frac{1}{\alpha} \displaystyle \sup_{v \in V_h \setminus \{0\}} 
\frac{|a_h(\tilde{u},v)-F_h(v)| + |a_h(\tilde{u}-w,v)|}{\parallel {\bf grad}\; v \parallel_{0,h}}.
\end{equation}
Thanks to the following facts the first term in the numerator of (\ref{fourthbound}) is expressed as in (\ref{ahFhprime}): Since $\tilde{u} \in H^3(\tilde{\Omega})$ 
we can apply First Green's identity to $a_h(\tilde{u},v)$ thereby getting rid of integrals on portions of $\Gamma$; next we note that 
$\nu \Delta u - {\bf b} \cdot {\bf grad}\; u + f=0$ in every $T \in {\mathcal T}_h \setminus {\mathcal S}_h$; this is also 
true of elements $T$ not belonging to the subset ${\mathcal Q}_h$ of ${\mathcal S}_h$ consisting of elements $T$ 
such that $T \setminus \Omega$ is not restricted to a set of vertexes of $\Omega_h$; finally we recall that $\nu \Delta \tilde{u} - 
{\bf b} \cdot {\bf grad}\; \tilde{u} + f$ vanishes identically in the set 
$T \cap \Omega$ and denote by $\tilde{\Delta}_T$ the interior of the set $T \setminus \Omega$ $\forall T \in {\mathcal Q}_h$. In short we can write:  
\begin{equation}
\label{ahFhprime}
|a_h(\tilde{u},v)-F_h(v)| = \displaystyle \sum_{T \in {\mathcal Q}_h} \left| \int_{\tilde{\Delta}_T} v(\nu \Delta \tilde{u} - {\bf b} \cdot {\bf grad} \; \tilde{u}) \right| 
\leq \displaystyle \sum_{T \in {\mathcal Q}_h} \parallel \nu \Delta \tilde{u} - {\bf b} \cdot {\bf grad}\; \tilde{u} \parallel_{0,\tilde{\Delta}_T} \parallel v \parallel_{0,\tilde{\Delta}_T} \!.
\end{equation}   
Now taking into account that $v \equiv 0$ on $\Gamma_h$ and recalling the constant $C_{\Gamma}$ defined in Lemma 2.1, it holds :   
$|v({\bf x})| \leq C_{\Gamma} h_T^2 \parallel |{\bf grad}\; v| \parallel_{0,\infty,\tilde{\Delta}_T}$, $\forall {\bf x} \in \tilde{\Delta}_T$, where 
$\parallel \cdot \parallel_{0,\infty,D}$ denotes the standard norm of $L^{\infty}(D)$, $D$ being a bounded open set of $\Re^2$. Now from a classical inverse  inequality we have $\parallel |{\bf grad} \; v| \parallel_{0,\infty,\tilde{\Delta}_T} \leq C_I h_T^{-1} \parallel {\bf grad} \; v \parallel_{0,T}$ for a mesh-independent constant $C_I$. Noticing that the measure of $\tilde{\Delta}_T$ is bounded by a constant depending only on $\Omega$ times $h_T^3$, after straightforward calculations we obtain for a certain mesh-independent constant $C_R$:
\begin{equation}
\label{fifthbound}
\parallel \nu \Delta \tilde{u} - {\bf b} \cdot {\bf grad}\; \tilde{u} \parallel_{0,\Delta^{'}_T} \parallel v \parallel_{0,\Delta^{'}_T} \leq C_R h_T^{5/2} \parallel \nu \Delta \tilde{u} - {\bf b} \cdot {\bf grad}\; \tilde{u} \parallel_{0,\tilde{\Delta}_T} \parallel 
{\bf grad} \; v \parallel_{0,T} \; \forall T \in {\mathcal Q}_h.
\end{equation}
Now plugging (\ref{fifthbound}) into (\ref{ahFhprime}) and applying the Cauchy-Schwarz inequality, we easily come up with,
\begin{equation}
\label{sixthbound}
|a_h(\tilde{u},v)-F_h(v)| \leq C_R h^{5/2} \displaystyle \parallel \nu \Delta \tilde{u} - {\bf b} \cdot {\bf grad}\; \tilde{u} \parallel_{0,\tilde{\Omega}} \parallel {\bf grad} \; v \parallel_{0,h}.
\end{equation}
Finally plugging (\ref{sixthbound}) into (\ref{fourthbound}) we immediately establish the validity of error estimate (\ref{estimateP2}). \QED \\ 

\begin{theorem} 
Let $k=3$ and assume that $u \in H^4(\Omega)$. Provided $h$ is sufficiently small, there exists a mesh independent constant $C_3$ such that the unique solution $u_h$ to (\ref{Poissonh}) satisfies:
\begin{equation}
\label{estimateP3} 
\parallel {\bf grad}(u - u_h) \parallel_{0,\tilde{\Omega}_h} \leq C_3 h^{3} [| \tilde{u} |_{4,\tilde{\Omega}} + \parallel \nu \Delta \tilde{u} - {\bf b} \cdot {\bf grad}\; \tilde{u} \parallel_{0,\infty,\tilde{\Omega}}]
\end{equation}
where $\tilde{u} \in H^4(\tilde{\Omega})$ is the regular extension of $u$ to $\tilde{\Omega}$ constructed in accordance to Stein et al. \cite{Stein}.
\end{theorem}

\prov
First of all we point out that, according to the Sobolev Embedding Theorem \cite{Adams}, $\Delta \tilde{u} \in L^{\infty}(\tilde{\Omega})$, since $\tilde{u} \in 
H^4(\tilde{\Omega})$ by assumption. \\ 
Now following the same steps as in the proof of Theorem \ref{P2} up to equation (\ref{ahFhprime}), the latter becomes for a certain mesh-independent constant $C_S$,
\begin{equation}
\label{ahFh3}
|a_h(\tilde{u},v)-F_h(v)| \leq C_S \displaystyle \sum_{T \in {\mathcal Q}_h}  h_T^{3} \parallel \nu \Delta \tilde{u} - {\bf b} \cdot {\bf grad}\; \tilde{u} \parallel_{0,\infty,\tilde{\Delta}_T} 
\parallel v \parallel_{0,\infty,\tilde{\Delta}_T}, 
\end{equation}
Akin to the previous proof, using a classical inverse inequality for triangles, we note that, 
\begin{equation}
\label{inverse} \parallel v \parallel_{0,\infty,\tilde{\Delta}_T} \leq C_{\Gamma} h_T^2 \parallel |{\bf grad}\; v| \parallel_{0,\infty,\tilde{\Delta}_T} 
\leq C_{\Gamma} C_I h_T \parallel {\bf grad}\; v \parallel_{0,T}. 
\end{equation}
Combining (\ref{ahFh3}) with (\ref{inverse}) we come up with,
\begin{equation}
\label{ahFh4}
|a_h(\tilde{u},v)-F_h(v)| \leq C_S C_{\Gamma} C_I \parallel \nu \Delta \tilde{u} - {\bf b} \cdot {\bf grad}\; \tilde{u} \parallel_{0,\infty,\tilde{\Omega}} \displaystyle \sum_{T \in {\mathcal Q}_h} h_T^4 
\parallel {\bf grad} \; v \parallel_{0,T},  
\end{equation}
Further applying the Cauchy-Schwarz inequality to the right hand side of (\ref{ahFh4}) we easily obtain:
\begin{equation}
\label{ahFh5}
|a_h(\tilde{u},v)-F_h(v)| \leq C_S C_{\Gamma} C_I h^3 \parallel\nu \Delta \tilde{u} - {\bf b} \cdot {\bf grad}\; \tilde{u} \parallel_{0,\infty,\tilde{\Omega}} 
\displaystyle \left[ \sum_{T \in {\mathcal Q}_h} h_T^{2} \right]^{1/2} 
\parallel {\bf grad} \; v  \parallel_{0,h}.
\end{equation}
From the assumptions on the mesh there exists a mesh-independent constant $C_J$ such that $[\sum_{T \in {\mathcal Q}_h} h_T^{2} ]^{1/2} $ $\leq C_J 
meas(\Gamma)$.  
Plugging this into (\ref{ahFh5}) and the resulting relation into (\ref{fourthbound}) we immediately establish error estimate (\ref{estimateP3}). \QED \\  

Akin to Theorem \ref{theorem1bis}, it is possible to establish error estimates in the $L^2$-norm in the case of a non convex $\Omega$, by requiring some more regularity from the solution $u$ of \eqref{Poisson}. However optimality is not attained, except for the case $k=2$. This is basically because of the absence of $u$ from the non-empty domain $\tilde{\Delta}_h : = \Omega_h \setminus \Omega$, whose area is an invariant $O(h^2)$ whatever $k$. Roughly speaking, integrals in $\tilde{\Delta}_h$ of expressions in terms of the approximate solution $u_h$ dominate the error, in such a way that the order of such terms cannot be reduced to less than $3.5$, even under additional regularity assumptions. \\
\indent Most steps in the proof of the following result rely on arguments essentially identical to those already exploited to prove Theorem \ref{theorem1bis}. Therefore we will focus on aspects specific to the non convex case.   
 
\begin{theorem}
\label{L2P2} Let $k=2$. Assume that $\Omega$ is not convex and $u \in H^{3+r}(\Omega)$ for 
$r=1/2+\epsilon$, $\epsilon>0 $ being arbitrarily small. Then provided $h$ is sufficiently small the following error estimate holds:
\begin{equation}
\label{L2estP2} \parallel u - u_h \parallel_{0,\tilde{\Omega}_h} \leq \tilde{C}_{0} h^{3} [G(\tilde{u})
+ \parallel u \parallel_{3+r,\Omega} ],   
\end{equation} 
where $\tilde{C}_{0}$ is a mesh-independent constant and $G(\tilde{u}):= | \tilde{u} |_{3,\tilde{\Omega}} 
+ h^{1/2} \parallel \nu \Delta \tilde{u} - {\bf b} \cdot {\bf grad} \; \tilde{u} \parallel_{0,\tilde{\Omega}}$.\\
\end{theorem}

\prov Let $\bar{u}_h$ be the function defined in $\Omega$ by $\bar{u}_h:=u_h-u$.\\ 
$v \in H^1_0(\Omega)$ being the function satisfying \eqref{adjoint}-\eqref{adjoint1}, we have: 
\begin{equation}
\label{L2estP2bis}
\parallel \bar{u}_h \parallel_{0,\tilde{\Omega}_h} \leq \parallel \bar{u}_h \parallel_{0} \leq 
C_{\Omega} \displaystyle \frac{- \int_{\Omega} \bar{u}_h(\nu \Delta v + {\bf b} \cdot {\bf grad} \; v)}{\parallel v \parallel_{2,\Omega}}.
\end{equation}
Now we set $\tilde{\Gamma}_h := \Omega_h \cap \Gamma$ and note that $meas(\tilde{\Gamma}_h) >0$. Using integration by parts we easily obtain,  
\begin{equation}
\label{L2estP2ter}
\parallel \bar{u}_h \parallel_{0,\tilde{\Omega}_h} \leq 
C_{\Omega} \displaystyle \frac{b_{1h}(\bar{u}_h,v) + \tilde{a}_h(\bar{u}_h,v)+a_{\Delta_h}(\bar{u}_h,v)}{\parallel v \parallel_{2,\Omega}},
\end{equation}
where the bilinear forms $b_{1h}$ and $a_{\Delta_h}$ are defined in \eqref{b1h} and \eqref{aDeltah}, respectively, and  
\begin{equation}
\label{tildeah}
\tilde{a}_h(w,z) 
:= \int_{\tilde{\Omega}_h} 
[\nu {\bf grad} \; w \cdot {\bf grad} \; z + ({\bf b} \cdot {\bf grad}\; w)z] \mbox{ for } w,z \in H^1(\Omega).
\end{equation}
On the other hand $\forall v_h \in V_h$ we have,
\begin{equation}
\label{L2estP2qua}
a_h(u_h,v_h)= \int_{\tilde{\Omega}_h} [-\nu \Delta u + {\bf b} \cdot {\bf grad} \;u]v_h =  
- \displaystyle \nu \int_{\tilde{\Gamma}_h} \frac{\partial u}{\partial n}v_h + \tilde{a}_h(u,v_h).
\end{equation}
Recalling the definitions of ${\mathcal Q}_h$ in the proof of Theorem \ref{P2} for every $T \in {\mathcal Q}_h$ and of the set $\tilde{\Delta}_T$ as the interior of $T \setminus \Omega$, we define 
\begin{equation}
\label{b5h}
b_{5h}(w,z) := \displaystyle \sum_{T \in {\mathcal Q}_h} \int_{\tilde{\Delta}_T} [\nu \Delta w - {\bf b} \cdot {\bf grad}\; w]z, 
\forall w \in W_h \mbox{ and } \forall z \in V_h.
\end{equation}
Denoting by $\tilde{\partial} T$ the set $\Gamma \cap T$ we further set,
\begin{equation}
\label{b6h}
b_{6h}(w,z) := \nu \displaystyle \sum_{T \in {\bf Q}_h} 
\int_{\tilde{\partial} T} \frac{\partial w}{\partial n}z \; \forall  
w \in W_h \cup H^2(\Omega) \mbox{ and } z \in V_h.
\end{equation}
It easily follows from \eqref{L2estP2qua} that 
\begin{equation}
\label{L2estP2qui}
-\tilde{a}_h(\bar{u}_h,v_h) + b_{5h}(u_h,v_h) + b_{6h}(\bar{u}_h,v_h)=0 \; \forall v_h \in V_h.
\end{equation}
Taking $v_h = \Pi_h(v)$, recalling that $e_h(v):=v-\Pi_h(v)$ and plugging \eqref{L2estP2qui} into \eqref{L2estP2ter} we come up with,
\begin{equation}
\label{L2estP2sex}
\parallel \bar{u}_h \parallel_{0,\tilde{\Omega}_h} \leq 
C_{\Omega} \displaystyle \frac{b_{1h}(\bar{u}_h,v) + 
+b_{5h}(u_h,\Pi_h(v))+b_{6h}(\bar{u}_h,\Pi_h(v))+\tilde{a}_h(\bar{u}_h,e_h(v))+a_{\Delta_h}(\bar{u}_h,v)}{\parallel v \parallel_{2,\Omega}}.
\end{equation}
On the other hand, recalling $b_{2h}$ given by \eqref{b2h} and using integration by parts we have
\begin{equation}
\label{L2estP2sept}
a_{\Delta_h}(\bar{u}_h,v)=a_{\Delta_h}(\bar{u}_h,e_h(v)) +  
 \displaystyle \sum_{T \in {\mathcal S}_h \setminus {\mathcal Q}_h} \int_{\partial T} \displaystyle 
\nu \frac{\partial \bar{u}_h}{\partial n} \Pi_h(v) + b_{2h}(\bar{u}_h,\Pi_h(v)).
\end{equation} 
Thus recalling $b_{3h}$ and $b_{4h}$ respectively defined by \eqref{b3h} and \eqref{b4h}, we finally obtain:
\begin{equation}
\label{L2estP2sex}
\left\{
\begin{array}{l}
\parallel \bar{u}_h \parallel_{0,\tilde{\Omega}_h} \leq 
C_{\Omega} \displaystyle \frac{L(\bar{u}_h,v)+b_{5h}(u_h,v_h)+\tilde{a}_h(\bar{u}_h,e_h(v))}{\parallel v \parallel_{2,\Omega}},\\
\mbox{where }\\
L(\bar{u}_h,v):=b_{1h}(\bar{u}_h,v) + b_{2h}(\bar{u}_h,\Pi_h(v)) + b_{3h}(\bar{u}_h,\Pi_h(v))+ 
b_{4h}(\bar{u}_h,e_h(v)).
\end{array}
\right.
\end{equation}
The estimation of $\tilde{a}_h(\bar{u}_h,e_h(v))$ is a trivial variant of the one in Theorem 
\ref{theorem1bis}, that is,
\begin{equation}
\label{estildeah}
\tilde{a}_h(\bar{u}_h,e_h(v)) \leq C_2 \tilde{C}_V h^{3} G(\tilde{u}) |v|_{2,\Omega},
\end{equation}
where $\tilde{C}_V$ is an interpolation error constant such that 
\begin{equation}
\label{intildeP1}
\parallel {\bf grad} [v-\Pi_h(v)]\parallel_{0,\tilde{\Omega}_h} \leq \tilde{C}_V h | v |_{2,\Omega}.
\end{equation}
\indent The bilinear forms $b_{ih}$, $i=1,2,3,4$ were studied in Theorem \ref{theorem1bis}. The corresponding estimates here are qualitatively the same taking $k=2$, if we replace here and there $| u |_{3,\Omega}$ 
by $G(u^{'})$. Hence all that is left to do is to estimate $b_{5h}(u_h,v_h)$. With this aim we proceed as follows:\\

Since $|v_h({\bf x})| \leq C_{\Gamma} h_T^2 \parallel {\bf grad}\; v_h \parallel_{0,\infty,T}$ $\forall {\bf x} \in \tilde{\Delta}_T$ and $\forall T \in {\mathcal Q}_h$ for $v_h \in V_h$, by a straightforward argument we can write
\begin{equation}
\label{step1}
b_{5h}(u_h,\Pi_h(v)) \leq \displaystyle \sum_{T \in {\mathcal Q}_h} C^2_{\Gamma} h_T^5 
[ \nu  \parallel \Delta u_h \parallel_{0,\infty,T} +
\parallel {\bf b} \parallel_{0,\infty} \parallel {\bf grad} \; u_h \parallel_{0,\infty,T} ]
\parallel {\bf grad}\; \Pi_h(v) \parallel_{0,\infty,T}.
\end{equation}
Since all the components of $[{\bf grad}\;\Pi_h(v)]_{|T}$ and $[H(u_h)]_{|T}$ are in ${\mathcal P}_0$ and those of $[{\bf grad}\;u_{h}]_{|T}$ are in ${\mathcal P}_1$, in all the norms involving $\Pi_h(v)$ and $u_h$ appearing in 
\eqref{step1} $T$ can be replaced by $T^{'}$. Thus by a classical inverse inequality and the Schwarz inequality we obtain successively, 
\begin{equation}
\label{step2}
b_{5h}(u_h,\Pi_h(v)) \leq \displaystyle \sum_{T \in {\mathcal Q}_h} C^2_{\Gamma} C_I^2 
 h_T^3 
[\nu \parallel \Delta u_h \parallel_{0,T^{'}} + \parallel {\bf b} \parallel_{0,\infty,\Omega}\parallel {\bf grad} \; u_h \parallel_{0,T^{'}} ]
\parallel {\bf grad}\; \Pi_h(v) \parallel_{0,T^{'}},
\end{equation}   
 
\begin{equation}
\label{step3}
b_{5h}(u_h,v_h) \leq  \sqrt{2} C^2_{\Gamma} C_I^2  h^3 \displaystyle \left\{ \sum_{T \in {\mathcal Q}_h}
\left[ \parallel \Delta u_h \parallel_{0,T^{'}}^2 + \parallel {\bf grad}\; u_h \parallel_{0,T^{'}}^2 \right]^{1/2} \right\} \parallel {\bf grad}\; \Pi_h(v) \parallel_{0,\tilde{\Omega}_h},
\end{equation}
On the other hand, by an inverse inequality in $T^{'}$ and owing to a classical approximation result, the Sobolev Embedding Theorem and an elementary geometric argument, there exists a mesh-independent constant $\bar{C}_I$ such that, 
\begin{equation}
\label{step4}
\begin{array}{l}
\parallel \Delta u_h \parallel_{0,T^{'}}^2 + \parallel {\bf grad} \; u_h \parallel_{0,T^{'}}^2 
\leq 3 \{\parallel \Delta(u_h-I_h(u)) \parallel_{0,T^{'}}^2 \\ 
+ \parallel {\bf grad}(u_h-I_h(u)) \parallel_{0,T^{'}}^2 + \parallel \Delta (I_h(u)-u)\parallel_{0,T^{'}}^2 +  \parallel {\bf grad}(I_h(u)-u)\parallel_{0,T^{'}}^2 +\parallel u \parallel_{2,T^{'}}^2 \} \\
\leq \bar{C}_I^2 
\{h_T^{-2} \parallel {\bf grad} [u_h-I_h(u)] \parallel_{0,T^{'}}^2 + h_T^2 | u |_{3,T^{'}}^2 + h_T^2 
\parallel u \parallel^2_{2,\infty,\Omega}\}.
\end{array}
\end{equation}
Plugging \eqref{step4} into \eqref{step3}, using the Cauchy-Schwarz inequality together with a simple trick for  the estimation of $b_{1h}$, taking into account 
\eqref{embedding} we easily obtain for another mesh-independent constant $\breve{C}_5$:
\begin{equation}
\label{step5}
b_{5h}(u_h,v_h) \leq  \breve{C}_5 h^3 \displaystyle \left\{ h G(\tilde{u}) + h^{1/2} \displaystyle \left[ \sum_{T \in {\mathcal Q}_h} h_T \right]^{1/2} \parallel u \parallel_{3+r,\Omega} \right\} \parallel {\bf grad}\; \Pi_h(v) \parallel_{0,\tilde{\Omega}_h}.
\end{equation}
On the other hand from \eqref{intildeP1} we easily derive, 
\begin{equation}
\label{interpol}
\parallel {\bf grad}\; \Pi_h(v) \parallel_{0,\tilde{\Omega}_h} \leq C_{\Pi} \parallel v \parallel_{2,\Omega}.
\end{equation} 
with $C_{\Pi}=\sqrt{1+\tilde{C}_V^2 diam(\Omega)^2}$. Hence there exists a mesh-independent constant 
$C_{b5}$ such that,
\begin{equation}
\label{estimateb5}
b_{5h}(u_h,v_h) \leq  C_{b5} h^{7/2} \displaystyle \left\{ h^{1/2} G(\tilde{u}) 
+ \parallel u \parallel_{3+r,\Omega} \right\} \parallel v  \parallel_{2,\Omega}.
\end{equation} 
Finally recalling \eqref{L2estP2sex} together with \eqref{estildeah}, \eqref{estimateb1}, \eqref{estimateb2}, \eqref{estimateb3} and \eqref{estimateb4}, estimate \eqref{estimateb5} completes the proof. \QED \\ 
              
\section{Numerical experiments}
In order to illustrate the error estimates derived in the previous section we solved equation (\ref{Poisson}) with our method in two test-cases, taking $k=2$. 
\subsection{Test-problem 1} 
\indent Here $\Omega$ is the ellipse delimited by the curve $(x/e)^2+y^2=1$ with $e>0$, $\nu=1$, ${\bf b} =(x,-y)$ and $d \equiv 0$, for an exact solution 
$u$ given by $u=(e^2-e^2 x^2 - y^2)(e^2-x^2-e^2 y^2)$. Thus we take $f := -\Delta u + {\bf b} \cdot {\bf grad} \; u$, and owing to symmetry we consider only the quarter domain given by $x>0$ and $y>0$ by prescribing Neumann boundary conditions on $x=0$ and $y=0$.  
We take $e=0.5$ and compute with quasi-uniform meshes defined by a single integer parameter $J$, constructed by a straightforward procedure. 
Roughly speaking the mesh of the quarter domain is the polar coordinate counterpart of the standard uniform mesh of the unit square $(0,1) \times (0,1)$ whose edges are parallel to the 
coordinate axes and to the line $x=y$.\\
In Table 1 we display the absolute errors in the norm $\parallel {\bf grad}(\cdot) \parallel_{0,h}$ and in the norm of $L^2(\Omega_h)$ 
for increasing values of $J$, more precisely $J=2^m$ for $m=2,3,4,5,6$. We also show the evolution of the maximum absolute errors at the mesh nodes denoted by $\parallel u - u_h \parallel_{0,\infty,h}$.\\  
\noindent As one infers from Table 1, the approximations obtained with our method perfectly conform to the theoretical estimate (\ref{errestconvex}). Indeed as $J$ increases the errors in the gradient $L^2$-norm decrease roughly as $(1/J)^2$, as predicted. The error in the $L^2$-norm in turn tends to decrease as $(1/J)^3$, while  
the maximum absolute error seem to behave like an $O(h^\beta)$, for $\beta$ less than but close to three.              
 
\begin{table*}[t!]
{\small 
\centering
\begin{tabular}{ccccccc} &\\ [-.3cm]  
$J$ & $\longrightarrow$ & 4 & $8$ & $16$ & $32$ & $64$ 
\tabularnewline & \\ [-.3cm] \hline &\\ [-.3cm]
$\parallel {\bf grad}(u-u_h) \parallel_{0,h}$ & $\longrightarrow$ & 0.539159 E-2 & 0.143611 E-2 & 0.367542 E-3 & 0.927845 E-4 & 0.233003 E-4   
\tabularnewline &\\ [-.3cm] \hline &\\ [-.3cm] 
$\parallel u-u_h \parallel_{0,h}$ & $\longrightarrow$ & 0.151255 E-3 & 0.184403 E-4 & 0.230467 E-5 & 0.289398 E-6 & 0.363189 E-7   
\tabularnewline &\\ [-.3cm] \hline &\\ [-.3cm]
$\parallel u-u_h \parallel_{0,\infty,h}$ & $\longrightarrow$ & 0.397339 E-3 & 0.751885 E-4 & 0.110067 E-4 & 0.148037 E-5 & 0.195523 E-6  
\tabularnewline & \\ [-.3cm] \hline &\\ [-.3cm]
\end{tabular}
\caption{Errors in different senses for Test-problem 1.}
}
\label{table1}
\end{table*}

\subsection{Test-problem 2}

\indent The aim of this Test-problem is to assess the behavior of our method in the case where $\Omega$ is non convex. Here we solve (\ref{Poisson}) for the following data: $\Omega$ is the annulus delimited by the circles given by $r=e < 1$ and $r=1$ with $r^2 = x^2+ y^2$, for an exact solution $u$ given by $\bar{u}=(r-e)(1-r)$ with $\bar{f}:=- \nu \Delta \bar{u}$, $\nu=1$, ${\bf b} \equiv {\bf 0}$ and $d \equiv 0$. Again we apply symmetry conditions on $x=0$ and $y=0$. 
We take $e=0.5$ and compute with quasi-uniform meshes defined by two integer parameters $I$ and $J$, constructed by subdividing the radial range $(0.5,1)$ into 
$J$ equal parts and the angular range $(0,\pi/2)$ into $I$ equal parts. In this way the mesh of the quarter domain is the polar coordinate counterpart of the 
$I \times J$ mesh of the rectangle $(0,\pi/2) \times (0.5,1)$ whose edges are parallel to the coordinate axes and to the line $x=\pi (y-0.5)$.\\
In Table 2 we display the absolute errors in the norm $\parallel {\bf grad}(\cdot) \parallel_{0,h}$ and in the norm of $L^2(\Omega_h)$ for $I=2J$,  
for increasing values of $I$, namely $I=2^m$ for $m=2,3,4,5,6$. We also show the evolution of the maximum absolute errors at the mesh nodes denoted by $\parallel u - u_h \parallel_{0,\infty,h}$.\\  
\noindent As one can observe, here again the quality of the approximations obtained with our method are in very good agreement with the theoretical result 
(\ref{errestconcave}), for as $J$ increases the errors in the gradient $L^2$-norm decrease roughly as $h^2$, as predicted. On the other hand here again the errors in the $L^2$-norm tend to decrease as $h^3$ and the maximum absolute errors behave like an $O(h^\beta)$ for $\beta$ close to and greater than three.      
           
\begin{table*}[t!]
{\small 
\centering
\begin{tabular}{ccccccc} &\\ [-.3cm]  
$I$ & $\longrightarrow$ & $4$ & $8$ & $16$ & $32$ & $64$ 
\tabularnewline & \\ [-.3cm] \hline &\\ [-.3cm]
$\parallel {\bf grad}(\bar{u}-u_h) \parallel_{0,h}$ & $\longrightarrow$ &  0.132906 E-1 & 0.334304 E-2 &  0.838061 E-3 & 0.209734 E-3 & 0.524545 E-4   
\tabularnewline &\\ [-.3cm] \hline &\\ [-.3cm] 
$\parallel \bar{u}-u_h \parallel_{0,h}$ & $\longrightarrow$ & 0.400090 E-3 & 0.491773 E-4 & 0.610753 E-5  & 0.761759 E-6 & 0.951819 E-7  
\tabularnewline &\\ [-.3cm] \hline &\\ [-.3cm]
$\parallel \bar{u}-u_h \parallel_{0,\infty,h}$ & $\longrightarrow$ & 0.679598E-3 & 0.716853 E-4 &  0.805631 E-5  &  0.947303 E-6 & 0.114561 E-6
\tabularnewline & \\ [-.3cm] \hline &\\ [-.3cm]
\end{tabular}
\caption{Errors in different senses for Test-problem 2.}
}
\label{table2}
\end{table*}

\section{Possible extensions and conclusions}

To conclude we make some comments on the methodology introduced in this work. We begin with general ones.  

\subsection{General considerations} 
\begin{enumerate}
\item
First of all it is important to stress that the assumption on the magnitude of the mesh parameter $h$ made throughout the paper is just a sufficient condition for the formal results given in this work to hold. It is by no means a necessary condition and actually we can even assert that it is rather an academic hypothesis. 
Indeed good numerical results can be obtained with meshes as coarse as can be. For example computations for test-problems like those given in Section 4 with the integer parameter $J=1$ or $J=2$ were carried out and no problem at all was detected.  
\item
The technique advocated in this work to solve the convection-diffusion equation in curved domains with classical Lagrange finite elements is actually much more general and universal. As a matter of fact it provides a simple and reliable manner to overcome technical difficulties brought about by more complicated problems and interpolations. This issue is illustrated in \cite{AMIS}, where we applied our technique to a Hermite analog of the Raviart-Thomas mixed finite element method of the lowest order to solve Maxwell's equations with Neumann boundary conditions. In a forthcoming paper we intend to complete this study 
by extending the technique to the Raviart-Thomas family \cite{RaviartThomas}, and to present the corresponding numerical analysis. As another example we can quote   
Hermite finite element methods to solve fourth order problems in curved domains with normal derivative degrees of freedom. Such d.o.f.s can also be dealt with very easily by means of our method, which is also shown in \cite{AMIS}. 
\item The solution of (\ref{Poisson}) with a non zero $d$ using our method is straightforward. Indeed, obviously enough, it suffices to assign the value of $d$ 
at each node belonging to the true boundary $\Gamma$ for any boundary element, that is, any element having an edge contained in $\Gamma_h$. The error estimates derived in this paper trivially extends to this case as the reader can certainly figure out. On the other hand in the case of Neumann boundary conditions $\partial u /\partial n = 0$  on $\Gamma$ (provided $f$ satisfies the underlying scalar condition) our method coincides with the standard 
Lagrange finite element method. Incidentally we recall that in case inhomogeneous Neumann boundary conditions are prescribed optimality can only be recovered if the linear form $F_h$ is modified in such a way that boundary integrals for boundary elements $T$ are shifted to the curved    
boundary portion of an element $\tilde{T}$ sufficiently close to the one of the corresponding curved element $T^{'}$. But this is an issue that has nothing to do 
with our method, which is basically aimed at resolving those related to the prescription of degrees of freedom for Dirichlet boundary conditions. 
\item As the reader has certainly noticed, in order to compute the element matrix and right hand side vector for a boundary element (in ${\mathcal S}_h$), we have to determine the inverse of an $n_k \times n_k$ matrix. However this extra effort should by no means be a problem at the current state-of-the art of Scientific Computing, as compared to the situation by the time isoparametric finite elements were introduced.
\item
It is important to recall that our method can do without numerical integration to compute element matrices, at least for quadratic and cubic finite elements, 
as pointed out in Sections 1 and 4. This is another significant advantage thereof over the isoparametric version of the finite element method. Indeed the latter  helplessly requires numerical integration for this purpose, since it deals with rational shape- and test-functions. While on the one hand this is not a real problem when the equation at hand is a simple one such as \eqref{Poisson}, on the other hand the choice of the right integration formula can turn to a sort of headache, in the case of more complex PDEs such as nonlinear ones.          
\end{enumerate}

\subsection{Comparison with the isoparametric technique}

The results in Section 4 validate the finite-element methodology studied in this article for the two-dimensional case, to solve boundary value problems posed in smooth curved domains. A priori it is an advantageous alternative in many respects to more classical techniques such as the isoparametric version of the finite element method. This is because its most outstanding features are not only universality but also simplicity, and eventually accuracy and CPU time too, although the two latter aspects were not our point from the beginning. 
Nevertheless we have compared our technique with the isoparametric one in both respects, by solving another test-problem using both approaches. It turned out that the new method was a little more accurate all the way. Just to illustrate this assertion we supply in Table 3 the errors in the $L^2(\Omega_h)$-norm of the solution gradient and of the solution itself, when both methods with $k=2$ are used to solve a toy Poisson equation $-\Delta u = f$ in the unit disk for 
$f(x,y):=9 (x^2+y^2)^{1/2}$ with $u = 0$ on the boundary. The exact solution is $u(x,y)=1-(x^2+y^2)^{3/2}$. The meshes employed in these computations are of the same type as those used in Test-problem 1 for an elliptical domain, i.e. they depend on an integer parameter $J$ in such a way that $h=1/J$. In Table 3 the solution obtained with isoparametric elements is denoted by $\tilde{u}^h$. Crout's method was employed for both methods to solve the resulting linear systems. \\
\indent Table 3 shows that the new method is a little more accurate than the isoparametric technique. In terms of CPU time the figures displayed in Table 4 are less conclusive. Indeed the new method can be considered globally less demanding than the isoparametic technique, though not uniformly (cf. the case $h=1/64$). As we should point out this comparison of CPU times is fair, since only boundary elements were treated differently for both methods, as required. \\
\begin{table*}[t!]
{\small 
\centering
\begin{tabular}{ccccccc} &\\ [-.3cm]  
$h$ & $\longrightarrow$ & 1/8 & $1/16$ & $1/32$ & $1/64$ & $1/128$ 
\tabularnewline & \\ [-.3cm] \hline &\\ [-.3cm]
$\parallel {\bf grad}_h(u-u_h) \parallel_{0,h}$ & $\longrightarrow$ & 0.361685 E-2 & 0.918504 E-3 & 0.231512 E-3 & 0.581281 E-4 & 0.145647 E-4    
\tabularnewline &\\ [-.3cm] \hline &\\ [-.3cm] 
$\parallel {\bf grad}(u-\tilde{u}_h) \parallel_{0,h}$ & $\longrightarrow$ &  0.383671 E-2 & 0.947667 E-3 & 0.235271 E-3 &  0.586053 E-4 & 0.146248 E-4 
\tabularnewline &\\ [-.3cm] \hline &\\ [-.3cm]
\tabularnewline &\\ [-.3cm] \hline &\\ [-.3cm]
$\parallel u-u_h \parallel_{0,h}$ & $\longrightarrow$ & 0.564603 E-4 & 0.717088 E-5 & 0.905923 E-6 & 0.124276 E-6 & 0.142626 E-7
\tabularnewline &\\ [-.3cm] \hline &\\ [-.3cm] 
$\parallel u-\tilde{u}_h \parallel_{0,h}$ & $\longrightarrow$ &  0.604713 E-4 & 0.744364 E-5 & 0.924795 E-6 &  0.128341 E-6 & 0.143329 E-7 
\tabularnewline &\\ [-.3cm] \hline &\\ [-.3cm]
\end{tabular} 
\caption{Errors with the new and the isoparametric approach for a test-problem in a disk taking $k=2$.} 
}
\label{table7}
\end{table*} 

\begin{table*}[t!]
{\small 
\centering
\begin{tabular}{ccccccc} &\\ [-.3cm]  
$h$ & $\longrightarrow$ & 1/8 & $1/16$ & $1/32$ & $1/64$ & $1/128$ 
\tabularnewline & \\ [-.3cm] \hline &\\ [-.3cm]
New approach & $\longrightarrow$ & 0.0384 & 0.2693 & 3.0017 & 61.8033 & 1117.1423  
\tabularnewline &\\ [-.3cm] \hline &\\ [-.3cm] 
Isoparametric approach & $\longrightarrow$ & 0.0367 &  0.3307 & 3.1674 & 46.4484 &  1201.2204   
\tabularnewline &\\ [-.3cm] \hline &\\ [-.3cm]
\end{tabular} 
\caption{CPU time in seconds to run a test-problem in a disk taking $k=2$.} 
}
\label{table8}
\end{table*}

\subsection{A short account of the three-dimensional case} 

\indent Saying a few words about the extremely important three-dimensional case is mandatory. \\

The three-dimensional counterpart of the method studied in this paper is studied in detail in \cite{arXiv3D}. 
Here we give only some highlights thereof.   
Although in this case too the method applies to much more general boundary value problems, for the sake of brevity we confined ourselves to the Poisson equation.\\ 
First of all for $N=3$ we make the very realistic assumption that an element $T \in {\mathcal T}_h$ has at most   
one face on $\Gamma_h$, and if no such a face exists $T$ has at most one edge on $\Gamma_h$. 
Actually we have to consider two subsets of ${\mathcal T}_h$, namely the subset ${\mathcal S}_h$ consisting of tetrahedra having one face on $\Gamma_h$ 
and the subset ${\mathcal R}_h$ consisting of tetrahedrons having exactly one edge on $\Gamma_h$. In contrast to the two-dimensional case, 
for $T \in {\mathcal S}_h$ it is not possible to define the set $\Delta_T$ delimited by $\Gamma$ and the face $F_T$ of $T$ contained in $\Gamma_h$, or equivalently the three skins associated with the three edges of $F_T$, in such a way that an underlying space $W_h$ of continuous functions is generated. Otherwise stated, 
in the three-dimensional case we have to deal with a non conforming space $W_h$. However this is not really a problem since the test-function space $V_h$ remains conforming. Nevertheless, at least from the formal point of view one had better employ a systematic way to extend or restrict the elements in ${\mathcal S}_h$ 
in order to construct a companion mesh of the whole $\Omega$ consisting of non overlapping straight elements $T \in {\mathcal T}_h \setminus {\mathcal S}_h$ and curved elements $T^{'}$ associated with $T \in {\mathcal S}_h$. Among other possibilities we can proceed as follows. For the latter elements, $T^{'}$ is delimited by $\Gamma$, the boundary portions of $T$ lying inside $\Omega$, and three skins $\delta_e$ corresponding to the three edges of the face $F_T \subset \Gamma_h$ generically denoted by $e$. $\delta_e$ lies on the plane containing $e$ that bisects the dihedral formed by two mesh faces whose intersection is $e$. Typically the pair of faces under consideration would correspond to the largest angle formed by two such faces. The interpolation points on $\Gamma$ pertaining to $T \in {\mathcal S}_h$ which are nodal points of $W_h$, are simply the intersections with $\Gamma$ of the perpendicular to $e$ in $\delta_e$ passing though the Lagrangian nodes 
of $e$. It is noteworthy that such nodes are interpolation nodes replacing Lagrangian nodes of $e$ for an element $T \in {\mathcal R}_h$ having $e$ as an edge,   although it is not necessary to consider any extension $T{'}$ of such a $T$. For every boundary mesh edge $e$ we denote by ${\mathcal L}_e$ the set of $k+1$ 
nodes belonging to $\overline{\delta}_e$ defined in the above manner.\\ 
This apparently complicated definition is aimed at ensuring that there is an extension ${\mathcal T}_h^{'}$ of the partition ${\mathcal T}_h$ consisting of non overlapping sets $T^{'}$ extending or restricting $T$, or doing both things at a time (typically $T^{'} := T \cup \Delta_T$ or $T^{'}: = \overline{T 
\setminus \Delta_T}$ according to the local configuration of $\Gamma$), besides the elements in ${\mathcal T}_h \setminus {\mathcal S}_h$. \\
Now for $w \in W_h$, $\forall T \in {\mathcal S}_h \cup {\mathcal R}_h$ and for every edge $e \subset T \cap \Gamma_h$, $w(P) = d(P)$ for all $P \in {\mathcal L}_e$. 
If $T \in {\mathcal R}_h$ all the remaining $(k+5)(k+1)k/6$ nodes used to define $w_{|T}$ for $w \in W_h$ are Lagrangian nodes of $T$. As for $T \in {\mathcal S}_h$, besides the $3k$ nodes in the three pertaining $\delta_e$s and its $(k+2)(k+1)k/6$ Lagrangian nodes  
not lying on $\Gamma_h$, for $k>2$ only, the remaining $(k-1)(k-2)/2$ nodes of $T \in {\mathcal S}_h$ associated with $W_h$ are the intersections with $\Gamma$ 
of the line passing through the vertex $O_T$ of $T$ not belonging to $\Gamma$ and the points  
subdividing the face opposite to $O_T$ into $k^2$ equal triangles, except those lying on the edges of $F_T$. Notice that, provided $h$ is small enough, there is no 
chance for two out of thus constructed $(k+3)(k+2)(k+1)/6$ nodes of $T \in {\mathcal S}_h \cup {\mathcal R}_h$ to be too close to each other, let alone to coincide. \\
Once the space $W_h$ is defined in accordance with the above constructions, the approximate problem (\ref{Poissonh}) can be posed in the same way as in the two-dimensional case. Corresponding existence, uniqueness and uniform stability results can be demonstrated in basically the same manner as in Section 2. As for error estimates, qualitative results equivalent to those proved in Section 3 can be expected to hold. Nonetheless their proof is at the price of several additional technicalities, especially in the non convex case. We address all those issues more thoroughly in \cite{arXiv3D}. 

\subsection{A word about further applications}
To close this work, we would like to insist that the technique advocated in this work to handle Dirichlet conditions prescribed on curvilinear boundaries has a wide scope of applicability. This is particularly true of some cases not so thoroughly addressed in the literature so far, such as mixed finite element methods for the incompressible Navier-Stokes equations. In this respect we refer to \cite{CAMWA}. Applications to elasticity problems can be found in \cite{CFM2017} and \cite{Maugin}.


\begin{thebibliography}{00}
\bibitem{Adams} 
R.A. Adams.
\newblock{\em Sobolev Spaces}.
Academic Press, N.Y., 1975.
\bibitem{Babuska}
I. Babu\v{s}ka.
\newblock{The finite element method with Lagrange multipliers}.
\newblock{\em Numerische Mathematik}, 20 (1973), 170--192.
\bibitem{RT}
F. Bertrand, S. M\"unzenmaier and G. Starke.
\newblock{First-order system least-squares on curved boundaries: higher-order Raviart--Thomas elements}.
\newblock{\em SIAM J. Numerical Analysis} 52-6 (2014), 3165-3180.
\bibitem{BrennerScott} 
S.C. Brenner and L.R.Scott. 
\newblock{\em The Mathematical Theory of Finite Element Methods}. Texts in Applied Mathematics 15, 
Springer, 2008.
\bibitem{Brezzi} F. Brezzi. 
\newblock{On the existence, uniqueness and approximation of saddle-point problems arising from 
Lagrange multipliers}. {\em RAIRO Analyse Num\'erique}. 8-2 (1974), 129-151.
\bibitem{BrezziFortin} 
F. Brezzi and M. Fortin (eds.). \newblock{\em Mixed and Hybrid Finite Element Methods}. 
Springer Series in Computational Mathematics, Vol. 15, 1991.
\bibitem{Ciarlet}
P.G. Ciarlet.
\newblock{\em The Finite Element Method for Elliptic Problems}.
North Holland, Amsterdam, 1978.
\bibitem{CiarletRaviart}
P.G. Ciarlet and P.A. Raviart. \newblock{The combined effect of curved boundaries and numerical integration in isoparametric finite element methods}. In: 
\newblock{\em The Mathematical Foundations of the Finite Element Method with Applications to Partial Differential Equations}, A.K. Aziz ed., pp. 409--474, 
Academic Press, 1972.
\bibitem{Coffman}
D. Coffman, D. Legg and Y. Pan. \newblock{A Taylor series condition for harmonic extension}. \newblock{\em  Real Analysis Exchange}, 28-1, (2002--2003), 235--253.
\bibitem{COAM}
J.A. Cuminato and V. Ruas.
{\em Unification of distance inequalities for linear variational problems}. 
{\em Computational and Applied Mathematics}, 34 (2015), 1009-1033. 
\newblock{\em Mathematics of Computation}, 77-261 (2008), 201--219.
\bibitem{Lions} 
J.-L. Lions. \newblock{\em Probl\`emes aux limites dans les \'equations aux d\'eriv\'ees partielles}, Presses de l'Universit\'e de Montr\'eal, Montr\'eal, 1962.
\bibitem{LionsMagenes} 
J.-L. Lions and E. Magen\`es. \newblock{\em Probl\`emes aux limites non homog\`enes et applications}, Dunod, Paris, 1968.
\bibitem{Necas}
J. Ne\v{c}as. \newblock{\em Les m\'ethodes directes en th\'eorie des \'equations elliptiques}, Masson, Paris, 1967.
\bibitem{Nitsche}
J. Nitsche. \newblock{On Dirichlet problems using subspaces with nearly zero boundary conditions}. 
\newblock{\em The Mathematical Foundations of the Finite Element Method with Applications to 
Partial Differential Equations}, A.K. Aziz ed., Academic Press, 1972.
\bibitem{Quarteroni}
A. Quarteroni, R. Sacco and F. Saleri, {\em Numerical Mathematics}, 
Texts in Applied Mathematics, Springer, 2007. 
\bibitem{RaviartThomas} 
P.-A.~Raviart and J.-M.~Thomas.
\newblock{Mixed Finite Element Methods for Second Order Elliptic Problems}. 
\newblock{\em Lecture Notes in Mathematics, Springer Verlag}, 606: 292--315, 1977.
\bibitem{arXiv3D}
V.~Ruas.  
\newblock{Methods of arbitrary optimal order with tetrahedral finite-element meshes forming polyhedral approximations of curved domains}. {\em arXiv Numerical Analysis}, arXiv:1706.08004 [math.NA], 2017.
\bibitem{CFM2017} V. Ruas. \newblock{A simple alternative for accurate finite-element modeling in curved domains}. \newblock{\em Comptes-rendus du Congr\`es Fran\c{c}ais de M\'ecanique}, Lille, France, 2017.
\bibitem{AMIS}
V. Ruas and M.A. Silva Ramos. 
\newblock{A Hermite Method for Maxwell's Equations}. \newblock{\em Applied Mathematics and Information Sciences}. 12-2 (2018), 271--283.
\bibitem{Maugin} 
V. Ruas. \newblock{Optimal Calculation of Solid-Body Deformations with Prescribed Degrees of Freedom over Smooth Boundaries}. In: \newblock{\em  
Advanced Structured Materials}, H. Altenbach, J. Pouget, M. Rousseau, B. Collet and T. Michelitsch (Org.), Magdeburg, 
Springer International Publishing, v.1, p. 695--704, 2018.
\bibitem{CAMWA} V. Ruas. \newblock{Accuracy enhancement for non-isoparametric finite-element simulations in curved domains; application to fluid flow}, \newblock{\em Computer and Mathematics with Applications}, to appear (link to on-line version https://doi.org/10.1016/j.camwa.2018.05.029).  
\bibitem{Scott}
L. R. Scott. \newblock{\em Finite Element Techniques for Curved Boundaries}. PhD thesis, MIT, 1973.
\bibitem{Stein}
D. B. Stein, R. D. Guy, B. Thomases. \newblock{Immersed boundary smooth extension:
A high-order method for solving PDE on arbitrary smooth
domains using Fourier spectral methods}. \newblock{\em Journal of Computational Physics}, 304 (2016), 252--274.
\bibitem{Strang}
G. Strang and G. Fix. \newblock{\em An Analysis of the Finite Element Method}. Prentice Hall, 1973.
\bibitem{Zenisek}
A. \v{Z}\'eni\v{s}ek. \newblock{Curved triangular finite $C^m$-elements}. \newblock{Aplikace Matematiky}, 23-5 (1978), 346--377.
\bibitem{Zienkiewicz} 
O.C. Zienkiewicz.
\newblock{\em The Finite Element Method in Engineering Science}. McGraw-Hill, 1971.
\bibitem{Zlamal} M. Zl\'amal. \newblock{Curved Elements in the Finite Element Method. I}. 
\newblock{\em SIAM Journal on Numerical Analysis}, 10-1 (1973), 229--240.

\end{thebibliography}
\end{document}